\documentclass[10pt, A4paper]{article}

\usepackage{graphicx,url}
\usepackage{amssymb,amsmath}
\usepackage{epstopdf}
\usepackage{amssymb,amsmath,amsthm}

\usepackage{setspace}

\usepackage{url,color}
\usepackage{lscape}

\numberwithin{equation}{section}

\begin{document}

\title{Poisson-type processes governed by fractional and higher-order
recursive differential equations}
\author{L.Beghin\thanks{%
Sapienza University of Rome} \and E.Orsingher\thanks{%
Sapienza University of Rome}}
\maketitle

\begin{abstract}
We consider some fractional extensions of the recursive differential
equation governing the Poisson process, i.e.%
\begin{equation*}
\frac{d}{dt}p_{k}(t)=-\lambda (p_{k}(t)-p_{k-1}(t)),\qquad k\geq 1,t>0
\end{equation*}%
by introducing fractional time-derivatives of order $\nu ,$ $2\nu ,...,n\nu $
. We show that the so-called \textquotedblleft Generalized Mittag-Leffler
functions\textquotedblright\ $E_{\alpha ,\beta }^{k}(x)$ (introduced by
Prabhakar [20]) arise as solutions of these equations. The corresponding
processes are proved to be renewal, with density of the intearrival times
(represented by Mittag-Leffler functions) possessing power, instead of
exponential, decay, for $t\rightarrow \infty .$ On the other hand, near the
origin the behavior of the law of the interarrival times drastically changes
for the parameter $\nu $ varying in $\left( 0,1\right] .$

For integer values of $\nu $, these models can be viewed as a higher-order
Poisson processes, connected with the standard case by simple and explict
relationships.


\end{abstract}

\section{Introduction}

Many well-known differential equations have been extended by introducing
fractional-order derivatives with respect to time (for instance, the heat,
wave and telegraph equations, as well as the higher-order heat-type
equations) or with respect to space (for instance, the equations involving
the Riesz fractional operator).

Fractional versions of the Poisson processes have been already presented and
studied in the literature: in [9] the so-called fractional master equation
was considered. A similar model was treated in [12], where the equation
governing the probability distribution of the homogeneous Poisson process
was modified, by introducing the Riemann-Liouville fractional derivative.
The results are given in analytical form, in terms of infinite series or
successive derivatives of Mittag-Leffler functions. We recall the definition
of the (two-parameters) Mittag-Leffler function:%
\begin{equation}
E_{\alpha ,\beta }(x)=\sum_{r=0}^{\infty }\frac{x^{r}}{\Gamma (\alpha
r+\beta )},\quad \alpha ,\beta \in \mathbb{C},\text{ }Re(\alpha ),%
Re(\beta )>0,\,x\in \mathbb{R},  \label{ml}
\end{equation}%
(see [20], §1.2).

A different definition of Poisson fractional process has been proposed by
Wang and Wen [30] and successively studied in [31]-[32]: in analogy to the
well-known fractional Brownian motion, the new process is defined as a
stochastic integral with respect to the Poisson measure. It displays
properties similar to fractional Brownian motion, such as self-similarity
(in the wide-sense) and long-range dependence.

Another approach was followed by Repin and Saichev [23]: they start by
generalizing, in a fractional sense, the distribution of the interarrival
times $U_{j}$ between two Poisson events. This is expressed, in terms of
Mittag-Leffler functions, as follows, for $\nu \in \left( 0,1\right] :$%
\begin{equation}
f(t)=\Pr \left\{ U_{j}\in dt\right\} /dt=-\frac{d}{dt}E_{\nu ,1}(-t^{\nu
})=\sum_{m=1}^{\infty }\frac{(-1)^{m+1}t^{\nu m-1}}{\Gamma \left( \nu
m\right) },\quad t>0  \label{ml2}
\end{equation}%
and coincides with the solution to the fractional equation%
\begin{equation}
\frac{d^{\nu }f(t)}{dt^{\nu }}=-f(t)+\delta (t),\quad t>0  \label{ml3}
\end{equation}%
where $\delta (\cdot )$ denotes the Dirac delta function and again the
fractional derivative is intended in the Riemann-Liouville sense. For $\nu
=1 $ formula (\ref{ml2}) reduces to the well-known density appearing in the
case of a homogeneous Poisson process, $N(t),t>0,$ with intensity $\lambda
=1 $, i.e. $f(t)=e^{-t}.$

The same approach is followed by Mainardi et al. [14]-[15]-[16], where a
deep analysis of the related process is performed: it turns out to be a true
renewal process, while it looses the Markovian property. Their first step is
the study of the following fractional equation (instead of (\ref{ml3}))%
\begin{equation}
\frac{d^{\nu }\psi (t)}{dt^{\nu }}=-\psi (t),  \label{ml4}
\end{equation}%
with initial condition $\psi (0^{+})=1$ and with fractional derivative
defined in the Caputo sense. The solution $\psi (t)=E_{\nu ,1}(-t^{\nu })$\
to (\ref{ml4}) represents the survival probability of the fractional Poisson
process. As a consequence its probability distribution is expressed in terms
of derivatives of Mittag-Leffler functions, while the density of the $k$-th
event waiting time is a fractional generalization of the Erlang distribution
and coincides with the $k$-fold convolution of (\ref{ml2}).

The analysis carried out by Beghin and Orsingher [2] starts, as in [11],
from the generalization of the equation governing the Poisson process, where
the time-derivative is substituted by the fractional derivative (in the
Caputo sense) of order $\nu \in \left( 0,1\right] $:%
\begin{equation}
\frac{d^{\nu }p_{k}}{dt^{\nu }}=-\lambda (p_{k}-p_{k-1}),\quad k\geq 0,
\label{ai}
\end{equation}%
with initial conditions%
\begin{equation*}
p_{k}(0)=\left\{
\begin{array}{c}
1\qquad k=0 \\
0\qquad k\geq 1%
\end{array}%
\right.
\end{equation*}%
and $p_{-1}(t)=0$. The main result is the expression of the solution as the
distribution of a composed process represented by the standard, homogeneous
Poisson process $N(t),t>0$ with a random time argument $\mathcal{T}_{2\nu
}(t),t>0$ as follows:%
\begin{equation*}
\mathcal{N}_{\nu }(t)=N(\mathcal{T}_{2\nu }(t)),\quad t>0.
\end{equation*}%
The process $\mathcal{T}_{2\nu }(t),t>0$ (independent of $N$) possesses a
well-known density, which coincides with the solution to a fractional
diffusion equation of order $2\nu $ (see (\ref{due.16}) below). In the
particular case where $\nu =1/2$ this equation coincides with the
heat-equation and the process representing time is the reflected Brownian
motion.

These results are reconsidered here, in the next section, from a different
point of view, which is based on the use of the Generalized Mittag-Leffler
(GML) function. The latter is defined as%
\begin{equation}
E_{\alpha ,\beta }^{\gamma }(z)=\sum_{r=0}^{\infty }\frac{\left( \gamma
\right) _{r}\,z^{r}}{r!\Gamma (\alpha r+\beta )},\quad \alpha ,\beta ,\gamma
\in \mathbb{C},\text{ }Re(\alpha ),Re(\beta ),Re(\gamma
)>0,  \label{gml2}
\end{equation}%
where $\left( \gamma \right) _{r}=\gamma (\gamma +1)...(\gamma +r-1)$ (for $%
r=1,2,...,$ and $\gamma \neq 0$) is the Pochammer symbol and $\left( \gamma
\right) _{0}=1$. The GML function has been extensively studied by Saxena et
al. (in [25]-[26]-[27] and [28]) and applied in connection with some
fractional diffusion equations, whose solutions are expressed as infinite
sums of (\ref{gml2}). For some properties of (\ref{gml2}), see also [29]. We
note that formula (\ref{gml2}) reduces to (\ref{ml}) for $\gamma =1.$

By using the function (\ref{gml2}) it is possible to write down in a more
compact form the solution to (\ref{ai}), as well as the density of the
waiting-time of the $k$-th event of the fractional Poisson process. As a
consequence some interesting relationships holding between the
Mittag-Leffler function (\ref{ml}) and the GML function (\ref{gml2}) are
obtained.

Moreover the use of GML functions allows us to derive an explicit expression
for the solution of the more complicated recursive differential equation,
where two fractional derivatives appear:%
\begin{equation}
\frac{d^{2\nu }p_{k}}{dt^{2\nu }}+2\lambda \frac{d^{\nu }p_{k}}{dt^{\nu }}%
=-\lambda ^{2}(p_{k}-p_{k-1}),\quad k\geq 0,  \label{due}
\end{equation}%
for $\nu \in \left( 0,1\right] $. As we will see in section 3, even in this
case we can define a process governed by (\ref{due}), which turns out to be
a renewal. The density of the interarrival times are no-longer expressed by
standard Mittag-Leffler functions as in the first case, but the use of GML
functions is required\ and the same is true for the $k$-th event
waiting-time.

An interesting relationship between the two models analyzed here can be
established by observing that the waiting-time of the $k$-th event of the
process governed by (\ref{due}) coincides in distribution with the waiting
time of the $(2k)$-th event for the first model. This suggests to interpret
our second model as a fractional Poisson process of the first type, which
jumps upward at even-order events $A_{2k}$ and the probability of the
successive odd-indexed events $A_{2k+1}$ is added to that of $A_{2k}$.
Correspondingly, the distribution of this second process $\widehat{\mathcal{N%
}}_{\nu }(t),t>0,$ can be expressed, in terms of the processes $N$ and $%
\mathcal{T}_{2\nu }$, as follows:%
\begin{equation*}
\Pr \left\{ \widehat{\mathcal{N}}_{\nu }(t)=k\right\} =\Pr \left\{ N(%
\mathcal{T}_{2\nu }(t))=2k\right\} +\Pr \left\{ N(\mathcal{T}_{2\nu
}(t))=2k+1\right\} ,\quad k\geq 0.
\end{equation*}%
We also study the probability generating functions of the two models, which
are themselves solutions to fractional equations; in particular in the
second case an interesting link with the fractional telegraph-type equation
is explored.

For $\nu =1$, equation (\ref{due}) takes the following form%
\begin{equation}
\frac{d^{2}p_{k}}{dt^{2}}+2\lambda \frac{dp_{k}}{dt}=-\lambda
^{2}(p_{k}-p_{k-1}),\quad k\geq 0
\end{equation}%
and the related process can be regarded as a standard Poisson process with
Gamma-distributed interarrival times (with parameters $\lambda ,2$). This is
tantamount to attribute the probability of odd-order values $A_{2k+1}$ of a
standard Poisson process to the events labelled by $2k.$ Moreover it should
be stressed that, in this special case, the equation satisfied by the
probability generating function $\widehat{G}(u,t),$ $t>0,$ $0<u\leq 1$, i.e.%
\begin{equation}
\frac{\partial ^{2}G(u,t)}{\partial t^{2}}+2\lambda \frac{\partial G(u,t)}{%
\partial t}=\lambda ^{2}(u-1)G(u,t),\qquad 0<\nu \leq 1
\end{equation}%
coincides with that of the damped oscillations.

All the previous results are further generalized to the case $n>2$ in the
concluding remarks: the structure of the process governed by the equation%
\begin{eqnarray}
&&\frac{d^{n\nu }p_{k}}{dt^{n\nu }}+\binom{n}{1}\lambda \frac{d^{(n-1)\nu
}p_{k}}{dt^{(n-1)\nu }}+...+\binom{n}{n-1}\lambda ^{n-1}\frac{d^{\nu }p_{k}}{%
dt^{\nu }}=-\lambda ^{n}(p_{k}-p_{k-1}),\quad k\geq 0,  \notag \\
&&\;  \label{nn}
\end{eqnarray}%
where $\nu \in \left( 0,1\right] $, is exactly the same as before and all
the previous considerations can be easily extended.

\section{First-type fractional recursive differential equation}

\subsection{The solution}

We begin by considering the following fractional recursive differential
equation
\begin{equation}
\frac{d^{\nu }p_{k}}{dt^{\nu }}=-\lambda (p_{k}-p_{k-1}),\quad k\geq 0,
\label{uno.2}
\end{equation}%
with $p_{-1}(t)=0$, subject to the initial conditions%
\begin{equation}
p_{k}(0)=\left\{
\begin{array}{c}
1\qquad k=0 \\
0\qquad k\geq 1%
\end{array}%
\right. .  \label{uno.1}
\end{equation}

We apply in (\ref{uno.2}) the definition of fractional derivative in the
sense of Caputo, that is, for $m\in \mathbb{N}$,%
\begin{equation}
\frac{d^{\nu }}{dt^{\nu }}u(t)=\left\{
\begin{array}{l}
\frac{1}{\Gamma (m-\nu )}\int_{0}^{t}\frac{1}{(t-s)^{1+\nu -m}}\frac{d^{m}}{%
ds^{m}}u(s)ds\text{,\qquad for }m-1<\nu <m \\
\frac{d^{m}}{dt^{m}}u(t)\text{,\qquad for }\nu =m%
\end{array}%
\right. .  \label{dc}
\end{equation}

We note that, for $\nu =1$, (\ref{uno.2}) coincides with the equation
governing the homogeneous Poisson process with intensity $\lambda >0$, so
that our first result generalizes the well-known distribution holding in the
standard case, i.e. $p_{k}(t)=\frac{(\lambda t)^{k}}{k!}e^{-\lambda t},$ $%
k\geq 0,t>0$.

We will obtain the solution to (\ref{uno.2})-(\ref{uno.1}) in terms of GML
functions (defined in (\ref{gml2})) and show that it represents a true
probability distribution of a process, which we will denote by $\mathcal{N}%
_{\nu }(t),t>0:$ therefore we will write%
\begin{equation}
p_{k}^{\nu }(t)=\Pr \left\{ \mathcal{N}_{\nu }(t)=k\right\} ,\quad \text{ }%
k\geq 0,t>0.  \label{bis}
\end{equation}

\

\noindent \textbf{Theorem 2.1 }The solution $p_{k}^{\nu }(t),$ for $%
k=0,1,... $\ and $t\geq 0,$\ of the Cauchy problem (\ref{uno.2})-(\ref{uno.1}%
) is given by%
\begin{equation}
p_{k}^{\nu }(t)=(\lambda t^{\nu })^{k}E_{\nu ,\nu k+1}^{k+1}(-\lambda t^{\nu
}),\quad \text{ }k\geq 0,t>0.  \label{gml3}
\end{equation}%
\textbf{Proof }We take the Laplace transform of equation (\ref{uno.2})
together with the condition (\ref{uno.1}) and consider that, for the Laplace
transform of the Caputo derivative, the following expression holds:%
\begin{eqnarray}
\mathcal{L}\left\{ \frac{d^{\nu }}{dt^{\nu }}u(t);s\right\}
&=&\int_{0}^{\infty }e^{-st}\frac{d^{\nu }}{dt^{\nu }}u(t)dt  \label{jen} \\
&=&s^{\nu }\mathcal{L}\left\{ u(t);s\right\} -\sum_{r=0}^{m-1}s^{\nu
-r-1}\left. \frac{d^{r}}{dt^{r}}u(t)\right\vert _{t=0},  \notag
\end{eqnarray}%
where $m=\left\lfloor \nu \right\rfloor +1$. Since, in this case, $\nu \in
\left( 0,1\right] ,$ we get $m=1.$ Therefore we get the following recursive
formula, for $k\geq 1,$%
\begin{equation}
(s^{\nu }+\lambda )\mathcal{L}\left\{ p_{k}^{\nu }(t);s\right\} =\lambda
\mathcal{L}\left\{ p_{k-1}^{\nu }(t);s\right\} ,  \label{rec}
\end{equation}%
while, for $k=0$, we obtain
\begin{equation}
\mathcal{L}\left\{ p_{0}^{\nu }(t);s\right\} =\frac{s^{\nu -1}}{s^{\nu
}+\lambda },
\end{equation}%
since equation (\ref{uno.2}) reduces to%
\begin{equation}
\frac{d^{\nu }p_{0}}{dt^{\nu }}=-\lambda p_{0},
\end{equation}%
with initial condition $p_{0}(0)=1.$

By applying (\ref{rec}) iteratively we get%
\begin{equation}
\mathcal{L}\left\{ p_{k}^{\nu }(t);s\right\} =\frac{\lambda ^{k}s^{\nu -1}}{%
(s^{\nu }+\lambda )^{k+1}}  \label{rec2}
\end{equation}%
which can be inverted by using formula (2.5) of [22], i.e.%
\begin{equation}
\mathcal{L}\left\{ t^{\gamma -1}E_{\beta ,\gamma }^{\delta }(\omega t^{\beta
});s\right\} =\frac{s^{\beta \delta -\gamma }}{(s^{\beta }-\omega )^{\delta }%
},  \label{pra}
\end{equation}%
(where $Re(\beta )>0,$ $Re(\gamma )>0,$ $Re(\delta )>0$
and $s>|\omega |^{\frac{1}{Re(\beta )}})$ for $\beta =\nu ,$ $\delta
=k+1$ and $\gamma =\nu k+1.$ Therefore the inverse of (\ref{rec2}) coincides
with (\ref{gml3}).\hfill $\square $

\

\noindent \textbf{Remark 2.1}\textit{\textbf{\ }As a cross check, we note
that, for }$\nu =1,$\textit{\ formula (\ref{gml3}) reduces to the
distribution of the homogeneous Poisson process since }$E_{1,k+1}^{k+1}(-%
\lambda t)=e^{-\lambda t}/k!.$

\textit{For any }$\nu \in \left( 0,1\right] $\textit{, it can be easily seen
that it coincides with the result, obtained by a different approach in [2], }%
\begin{eqnarray}
&&p_{k}^{\nu }(t)=(\lambda t^{\nu })^{k}\sum_{r=0}^{\infty }\frac{\left(
k+1\right) _{r}(-\lambda t^{\nu })^{r}}{r!\Gamma (\nu r+\nu k+1)}
\label{rec3} \\
&=&\frac{(\lambda t^{\nu })^{k}}{k!}\sum_{r=0}^{\infty }\frac{(r+k)!}{%
r!\Gamma \left( \nu (k+r)+1\right) }(-\lambda t^{\nu })^{r}  \notag \\
&=&\sum_{r=k}^{\infty }\binom{r}{k}\frac{(-1)^{r-k}(\lambda t^{\nu })^{r}}{%
\Gamma \left( \nu r+1\right) },\quad \text{ }k\geq 0,t>0,  \notag
\end{eqnarray}%
\textit{by noting that }%
\begin{eqnarray*}
(r+k)! &=&(k+1+r-1)(k+1+r-2)\cdot \cdot \cdot (k+1)k! \\
&=&\left( k+1\right) _{r}\,k!.
\end{eqnarray*}

\

\noindent \textbf{Remark 2.2}\textit{\ The result of Theorem 2.1 shows that
the first model proposed by Mainardi et al. [14] as a fractional version of
the Poisson process (called renewal process of Mittag-Leffler type)
coincides with the solution of equation (\ref{uno.2}) and therefore with (%
\ref{gml3}). In the paper cited above the distribution is expressed in terms
of successive derivatives of Mittag-Leffler functions as }%
\begin{equation}
v_{k}^{\nu }(t)=\frac{t^{k\nu }}{k!}\left[ \frac{d^{k}}{dx^{k}}E_{\nu
,1}\left( x\right) \right] _{x=-t^{\nu }},\quad k\geq 0,t>0  \label{mai}
\end{equation}%
\textit{and is obtained by means of the fractional generalization of the
Erlang density}%
\begin{equation}
f_{k}^{\nu }(t)=\frac{\nu t^{k\nu -1}}{(k-1)!}\left[ \frac{d^{k}}{dx^{k}}%
E_{\nu ,1}\left( x\right) \right] _{x=-t^{\nu }},  \label{maif}
\end{equation}%
\textit{which represents the distribution of the waiting-time of the }$k$%
\textit{-th event. Clearly, for }$\nu =1$\textit{, }$f_{k}(t)=\frac{%
t^{k-1}e^{-t}}{(k-1)!}$\textit{\ is the Erlang density representing the
waiting-time of the }$k$\textit{-th event for the standard Poisson process
(for }$\lambda =1$\textit{). Since we can rewrite (\ref{mai}) as}%
\begin{eqnarray*}
v_{k}^{\nu }(t) &=&\frac{t^{k\nu }}{k!}\left[ \frac{d^{k}}{dx^{k}}%
\sum_{j=0}^{\infty }\frac{x^{j}}{\Gamma (\nu j+1)}\right] _{x=-t^{\nu }} \\
&=&\frac{t^{k\nu }}{k!}\sum_{j=k}^{\infty }\frac{j(j-1)...(j-k+1)(-t^{\nu
})^{j-k}}{\Gamma (\nu j+1)} \\
&=&\frac{t^{k\nu }}{k!}\sum_{l=0}^{\infty }\frac{(l+k)!(-t^{\nu })^{l}}{%
l!\Gamma (\nu l+\nu k+1)},
\end{eqnarray*}%
\textit{it is evident that it coincides with (\ref{gml3}), for }$\lambda =1.$

\

We derive now an interesting relationship between the GML function in (\ref%
{gml3}) and the solution to a fractional-diffusion equation, which is
expressed in terms of the Wright function%
\begin{equation*}
W_{\alpha ,\beta }(x)=\sum_{k=0}^{\infty }\frac{x^{k}}{k!\Gamma (\alpha
k+\beta )},\qquad \alpha >-1,\text{ }\beta >0,\text{ }x\in \mathbb{R}.
\end{equation*}

Let us denote by $v_{2\nu }=v_{2\nu }(y,t)$ the solution to the Cauchy
problem
\begin{equation}
\left\{
\begin{array}{l}
\frac{\partial ^{2\nu }v}{\partial t^{2\nu }}=\lambda ^{2}\frac{\partial
^{2}v}{\partial y^{2}},\qquad t>0,\text{ }y\in \mathbb{R} \\
v(y,0)=\delta (y),\qquad \text{for }0<\nu <1 \\
v_{t}(y,0)=0,\qquad \text{for }1/2<\nu <1%
\end{array}%
\right. .  \label{due.16}
\end{equation}%
then it is well-known (see [17], p.142) that the solution of (\ref{due.16})
can be written as%
\begin{equation}
v_{2\nu }(y,t)=\frac{1}{2\lambda t^{\nu }}W_{-\nu ,1-\nu }\left( -\frac{|y|}{%
\lambda t^{\nu }}\right) ,\qquad t>0,\text{ }y\in \mathbb{R}.
\label{due.16b}
\end{equation}%
If we fold the above solution and define%
\begin{equation}
\overline{v}_{2\nu }(y,t)=\left\{
\begin{array}{l}
2v_{2\nu }(y,t),\qquad \text{ }y>0 \\
0,\qquad y<0%
\end{array}%
\right.  \label{due.16f}
\end{equation}%
then, for $k\geq 0,$ we get%
\begin{eqnarray}
p_{k}^{\nu }(t) &=&\int_{0}^{+\infty }e^{-y}\frac{y^{k}}{k!}\overline{v}%
_{2\nu }(y,t)dy  \label{due.19} \\
&=&\Pr \left\{ N(\mathcal{T}_{2\nu }(t))=k\right\} .  \notag
\end{eqnarray}%
In (\ref{due.19}) $T_{2\nu }(t),t>0$\ represents a random time with
transition density given in (\ref{due.16b})-(\ref{due.16f}) and independent
of the Poisson process $N(t),t>0$ (note that from now on $N$ denotes a
standard Poisson process with intensity $\lambda =1$). This result is
analogous to what happens for other kinds of fractional equations, such as
the diffusion ones (see [18]), the telegraph-type fractional equations (in
[17]) and even the higher-order heat-type equations with fractional
time-derivative (see [1]).

Formula (\ref{due.19}) was obtained for the first time in [2] and we prove
it here in an alternative form, by resorting to the Laplace transform. We
compare (\ref{rec2}) with the Laplace transform of the distribution $%
p_{k}(t),t>0$ of a standard Poisson process with intensity $\lambda >0$,
which reads%
\begin{equation}
\mathcal{L}\left\{ p_{k}(t);s\right\} =\frac{\lambda ^{k}}{(s+\lambda )^{k+1}%
}.  \label{rec7}
\end{equation}%
Formula (\ref{rec2}) can be consequently written as%
\begin{equation}
\mathcal{L}\left\{ p_{k}^{\nu }(t);s\right\} =s^{\nu -1}\mathcal{L}\left\{
p_{k}(t);s^{\nu }\right\} ,
\end{equation}%
which, by inversion, leads to the following convolution:%
\begin{equation}
p_{k}^{\nu }(t)=\frac{1}{\Gamma (1-\nu )}\int_{0}^{t}(t-w)^{-\nu }\mathcal{L}%
^{-1}\left\{ \int_{0}^{+\infty }e^{-s^{\nu }y/\lambda }\overline{p}%
_{k}(y)dy;w\right\} dw,  \label{pk}
\end{equation}%
where $\overline{p}_{k}(t),t>0$ represents the distribution of the Poisson
process with intensity $\lambda =1.$

The inverse transform in (\ref{pk}) can be expressed as follows%
\begin{equation}
\mathcal{L}^{-1}\left\{ \int_{0}^{+\infty }e^{-s^{\nu }y/\lambda }\overline{p%
}_{k}(y)dy;w\right\} =\int_{0}^{+\infty }g_{\nu }(w;\frac{y}{\lambda })%
\overline{p}_{k}(y)dy,  \label{pk1}
\end{equation}%
by recalling that
\begin{equation}
e^{-s^{\nu }y/\lambda }=\int_{0}^{+\infty }e^{-sz}g_{\nu }(z;\frac{y}{%
\lambda })dz,\qquad 0<\nu <1,\;w>0,  \label{es}
\end{equation}%
where $g_{\nu }(\cdot ;y)$ is a stable law $S_{\nu }(\mu ,\beta ,\sigma )$
of order $\nu $, with parameters $\mu =0,$ $\beta =1$ and $\sigma =\left(
\frac{y}{\lambda }\cos \frac{\pi \nu }{2}\right) ^{\frac{1}{\nu }}.$

By inserting (\ref{pk1}) into (\ref{pk}) we get
\begin{eqnarray}
p_{k}^{\nu }(t) &=&\frac{1}{\Gamma (1-\nu )}\int_{0}^{t}(t-w)^{-\nu }\left(
\int_{0}^{+\infty }g_{\nu }(w;\frac{y}{\lambda })\overline{p}%
_{k}(y)dy\right) dw  \label{pk2} \\
&=&\int_{0}^{+\infty }\left( \frac{1}{\Gamma (1-\nu )}\int_{0}^{t}(t-w)^{-%
\nu }g_{\nu }(w;\frac{y}{\lambda })dw\right) \overline{p}_{k}(y)dy.  \notag
\end{eqnarray}%
We recognize in (\ref{pk2}) the fractional integral of order $\nu $ of the
stable law $g_{\nu }$ and we apply the result (3.5) of [17], which can be
rewritten, in this case, as%
\begin{equation*}
\frac{1}{\Gamma (1-\nu )}\int_{0}^{t}(t-w)^{-\nu }g_{\nu }(w;\frac{y}{%
\lambda })dw=\overline{v}_{2\nu }(y,t),
\end{equation*}%
where $\overline{v}_{2\nu }(y,t)$ is given in (\ref{due.16f}).

We can conclude from (\ref{gml3}) that%
\begin{equation}
E_{\nu ,\nu k+1}^{k+1}(-\lambda t^{\nu })=\frac{1}{k!\lambda ^{k+1}t^{\nu
(k+1)}}\int_{0}^{+\infty }e^{-y}y^{k}W_{-\nu ,1-\nu }(-\frac{y}{\lambda
t^{\nu }})dy.  \label{gml4}
\end{equation}

\

\noindent \textbf{Remark 2.3 }\textit{Formula (\ref{due.19}) can be useful
also in checking that the sum of }$p_{k}^{\nu }(t)$\textit{\ (either in the
form (\ref{gml3}) or (\ref{rec3})), for }$k\geq 0$\textit{, is equal to one,
since it is }%
\begin{equation}
\sum_{k=0}^{\infty }p_{k}^{\nu }(t)=\int_{0}^{+\infty }e^{-y}\left(
\sum_{k=0}^{\infty }\frac{y^{k}}{k!}\right) \overline{v}_{2\nu
}(y,t)dy=\int_{0}^{+\infty }\overline{v}_{2\nu }(y,t)dy=1.  \notag
\end{equation}%
\textit{Moreover, since result (\ref{gml4}) holds for any }$t>0$\textit{, we
can choose }$t=1$\textit{, so that we get, by a change of variable,}%
\begin{equation*}
E_{\nu ,\nu k+1}^{k+1}(-\lambda )=\frac{1}{k!}\int_{0}^{+\infty }e^{-\lambda
y}y^{k}W_{-\nu ,1-\nu }(-y)dy.
\end{equation*}%
\textit{This shows that the GML function }$E_{\nu ,\nu k+1}^{k+1}$\textit{\
can be interpreted as the Laplace transform of the function }$\frac{y^{k}}{k!%
}W_{-\nu ,1-\nu }(-y).$

\textit{In particular, for }$\nu =\frac{1}{2},$\textit{\ since (\ref{due.19}%
) reduces to}

\begin{eqnarray}
\Pr \left\{ \mathcal{N}_{1/2}(t)=k\right\} &=&\int_{0}^{+\infty }e^{-y}\frac{%
y^{k}}{k!}\frac{e^{-y^{2}/4\lambda ^{2}t}}{\sqrt{\pi \lambda ^{2}t}}dy
\label{due.20} \\
&=&\Pr \left\{ N(|B_{\lambda }(t)|)=k\right\} ,  \notag
\end{eqnarray}%
\textit{where }$B_{\lambda }(t)$\textit{\ is a Brownian motion with variance
}$2\lambda ^{2}t$\textit{\ (independent of }$N)$, \textit{we get (for }$t=1$%
\textit{)}%
\begin{equation}
E_{\frac{1}{2},\frac{k}{2}+1}^{k+1}(-\lambda )=\frac{1}{k!}\int_{0}^{+\infty
}e^{-\lambda y}y^{k}\frac{e^{-y^{2}/4}}{\sqrt{\pi }}dy.  \label{gml5}
\end{equation}%
\textit{The previous relationship can be checked directly, as follows}%
\begin{eqnarray*}
\frac{1}{k!}\int_{0}^{+\infty }e^{-\lambda y}y^{k}\frac{e^{-y^{2}/4}}{\sqrt{%
\pi }}dy &=&\frac{1}{\sqrt{\pi }k!}\sum_{r=0}^{\infty }\frac{\left( -\lambda
\right) ^{r}}{r!}\int_{0}^{+\infty }y^{k+r}e^{-y^{2}/4}dy \\
&=&\frac{1}{\sqrt{\pi }k!}\sum_{r=0}^{\infty }\frac{\left( -\lambda \right)
^{r}2^{k+r}}{r!}\Gamma \left( \frac{r+k}{2}+\frac{1}{2}\right) \\
&=&\left[ \text{by the duplication formula}\right] \\
&=&\frac{1}{k!}\sum_{r=0}^{\infty }\frac{\left( -\lambda \right) ^{r}}{r!}%
\frac{\Gamma \left( r+k\right) (r+k)}{\frac{r+k}{2}\Gamma \left( \frac{r+k}{2%
}\right) } \\
&=&\frac{1}{k!}\sum_{r=0}^{\infty }\frac{\left( -\lambda \right) ^{r}}{r!}%
\frac{\left( r+k\right) !}{\Gamma \left( \frac{r+k}{2}+1\right) }=E_{\frac{1%
}{2},\frac{k}{2}+1}^{k+1}(-\lambda ).
\end{eqnarray*}

\subsection{Properties of the corresponding process}

From the previous results we can conclude that the GML function $E_{\nu ,\nu
k+1}^{k+1}(-\lambda t^{\nu })$, $k\geq 0,$ suitably normalized by the factor
$(\lambda t^{\nu })^{k}$, represents a proper probability distribution and
we can indicate it as $\Pr \left\{ \mathcal{N}_{\nu }(t)=k\right\} .$

Moreover by (\ref{due.19}) we can consider the process $\mathcal{N}_{\nu
}(t),t>0$ as a time-changed Poisson process. It is well-known (see [8] and
[11]) that, for a homogeneous Poisson process $N$ subject to a random time
change (by the random function $\Lambda (\left( 0,t\right] )$), the
following equality in distribution holds:%
\begin{equation}
N(\Lambda (\left( 0,t\right] ))\overset{d}{=}M(t),  \label{gra}
\end{equation}%
where $M(t),t>0$ is a Cox process directed by $\Lambda .$ In our case the
random measure $\Lambda (\left( 0,t\right] )$ possesses distribution $%
\overline{v}_{2\nu }$ given in (\ref{due.16b})-(\ref{due.16f}) and we can
conclude that $\mathcal{N}_{\nu }$ is a Cox process. This conclusion will be
confirmed by the analysis of the factorial moments.

Moreover, as remarked in [2] and [15], the fractional Poisson process $%
\mathcal{N}_{\nu }(t),t>0$ represents a renewal process with interarrival
times $U_{j}$ distributed according to the following density, for $j=1,2,...$%
:%
\begin{equation}
f_{1}^{\nu }(t)=\Pr \left\{ \mathcal{U}_{j}\in dt\right\} /dt=\lambda t^{\nu
-1}E_{\nu ,\nu }(-\lambda t^{\nu }),  \label{maimai}
\end{equation}%
with Laplace transform
\begin{equation}
\mathcal{L}\left\{ f_{1}^{\nu }(t);s\right\} =\frac{\lambda }{s^{\nu
}+\lambda }.  \label{mailap}
\end{equation}%
Therefore the density of the waiting time of the $k$-th event, $%
T_{k}=\sum_{j=1}^{k}U_{j}$, possesses the following Laplace transform%
\begin{equation}
\mathcal{L}\left\{ f_{k}^{\nu }(t);s\right\} =\frac{\lambda ^{k}}{(s^{\nu
}+\lambda )^{k}}.  \label{fk2}
\end{equation}%
Its inverse can be obtained by applying again (\ref{pra}) for $\beta =\nu ,$
$\gamma =\nu k$ and $\omega =-\lambda $ and can be expressed, as for the
probability distribution, in terms of a GML function as%
\begin{equation}
f_{k}^{\nu }(t)=\Pr \left\{ T_{k}\in dt\right\} /dt=\lambda ^{k}t^{\nu
k-1}E_{\nu ,\nu k}^{k}(-\lambda t^{\nu }).  \label{fk}
\end{equation}

Formula (\ref{fk}) coincides with (\ref{maif}), for $\lambda =1.$ The
corresponding distribution function can be obtained in two different ways.
The first one is based on (\ref{fk}) and yields%
\begin{eqnarray}
F_{k}^{\nu }(t) &=&\Pr \left\{ T_{k}<t\right\}  \label{fk5} \\
&=&\lambda ^{k}\int_{0}^{t}s^{\nu k-1}\sum_{j=0}^{\infty }\frac{%
(k-1+j)!(-\lambda s^{\nu })^{j}}{j!(k-1)!\Gamma (\nu j+\nu k)}ds  \notag \\
&=&\frac{\lambda ^{k}t^{\nu k}}{\nu }\sum_{j=0}^{\infty }\frac{%
(k-1+j)!(-\lambda t^{\nu })^{j}}{j!(k-1)!(k+j)\Gamma (\nu j+\nu k)}  \notag
\\
&=&\lambda ^{k}t^{\nu k}\sum_{j=0}^{\infty }\frac{(k-1+j)!(-\lambda t^{\nu
})^{j}}{j!(k-1)!\Gamma (\nu j+\nu k+1)}  \notag \\
&=&\lambda ^{k}t^{\nu k}E_{\nu ,\nu k+1}^{k}(-\lambda t^{\nu }).  \notag
\end{eqnarray}%
We can check that (\ref{fk5}) satisfies the following relationship%
\begin{equation}
\Pr \left\{ T_{k}<t\right\} -\Pr \left\{ T_{k+1}<t\right\} =p_{k}^{\nu }(t).
\label{fk8}
\end{equation}%
Indeed from (\ref{fk5}) we can rewrite (\ref{fk8}) as%
\begin{eqnarray*}
&&\lambda ^{k}t^{\nu k}E_{\nu ,\nu k+1}^{k}(-\lambda t^{\nu })-\lambda
^{k+1}t^{\nu (k+1)}E_{\nu ,\nu (k+1)+1}^{k+1}(-\lambda t^{\nu }) \\
&=&\lambda ^{k}t^{\nu k}\sum_{j=0}^{\infty }\frac{(k-1+j)!(-\lambda t^{\nu
})^{j}}{j!(k-1)!\Gamma (\nu j+\nu k+1)}-\lambda ^{k+1}t^{\nu
(k+1)}\sum_{j=0}^{\infty }\frac{(k+j)!(-\lambda t^{\nu })^{j}}{j!k!\Gamma
(\nu j+\nu k+\nu +1)} \\
&=&\left[ \text{by putting }l=j+1\text{ in the second sum}\right] \\
&=&\lambda ^{k}t^{\nu k}\sum_{j=0}^{\infty }\frac{(k-1+j)!(-\lambda t^{\nu
})^{j}}{j!(k-1)!\Gamma (\nu j+\nu k+1)}+\lambda ^{k}t^{\nu
k}\sum_{l=1}^{\infty }\frac{(k+l-1)!(-\lambda t^{\nu })^{l}}{(l-1)!k!\Gamma
(\nu l+\nu k+1)} \\
&=&\lambda ^{k}t^{\nu k}\sum_{j=0}^{\infty }\frac{(k+j)!(-\lambda t^{\nu
})^{j}}{j!k!\Gamma (\nu j+\nu k+1)}=p_{k}^{\nu }(t).
\end{eqnarray*}

The second method of evaluating the distribution function resorts to the
probabilities given in the form (\ref{rec3}):
\begin{eqnarray}
\Pr \left\{ T_{k}<t\right\} &=&\sum_{m=k}^{\infty }p_{m}^{\nu }(t)
\label{fk7} \\
&=&\sum_{m=k}^{\infty }\sum_{r=m}^{\infty }\binom{r}{m}\frac{%
(-1)^{r-m}(\lambda t^{\nu })^{r}}{\Gamma \left( \nu r+1\right) }  \notag \\
&=&\sum_{r=k}^{\infty }\frac{(-1)^{r}(\lambda t^{\nu })^{r}}{\Gamma \left(
\nu r+1\right) }\sum_{m=k}^{r}\binom{r}{m}(-1)^{m}.  \notag
\end{eqnarray}%
Finally if we rewrite (\ref{fk5}) as%
\begin{equation*}
\lambda ^{k}t^{\nu k}\sum_{j=k}^{\infty }\binom{j-1}{k-1}\frac{%
(-1)^{k}(-\lambda t^{\nu })^{j}}{\Gamma (\nu j+1)}
\end{equation*}%
and compare it with (\ref{fk7}), we extract the following useful
combinatorial relationship:%
\begin{equation*}
\binom{j-1}{k-1}(-1)^{k}=\sum_{m=k}^{j}\binom{j}{m}(-1)^{m},\qquad j\geq k.
\end{equation*}

\

\noindent \textbf{Remark 2.4 }\textit{As pointed out in [5] and [23], the
density of the interarrival times (\ref{maimai}) possess the following
asymptotic behavior, for }$t\rightarrow \infty $\textit{:}%
\begin{eqnarray}
\Pr \left\{ \mathcal{U}_{j}\in dt\right\} /dt &=&\lambda t^{\nu -1}E_{\nu
,\nu }(-\lambda t^{\nu })=-\frac{d}{dt}E_{\nu ,1}(-\lambda t^{\nu })
\label{rai} \\
&=&\lambda ^{1/\nu }\frac{\sin \left( \nu \pi \right) }{\pi }%
\int_{0}^{+\infty }\frac{r^{\nu }e^{-\lambda ^{1/\nu }rt}}{r^{2\nu }+2r^{\nu
}\cos (\nu \pi )+1}dr  \notag \\
&\sim &\frac{\sin \left( \nu \pi \right) }{\pi }\frac{\Gamma (\nu +1)}{%
\lambda t^{\nu +1}}=\frac{\nu }{\lambda \Gamma (1-\nu )t^{\nu +1}},  \notag
\end{eqnarray}%
\textit{where the well-known expansion of the Mittag-Leffler function}%
\begin{equation}
E_{\nu ,1}(-\lambda t^{\nu })=\frac{\sin \left( \nu \pi \right) }{\pi }%
\int_{0}^{+\infty }\frac{r^{\nu -1}e^{-\lambda ^{1/\nu }rt}}{r^{2\nu
}+2r^{\nu }\cos (\nu \pi )+1}dr  \label{rai3}
\end{equation}%
\textit{has been applied (see the Appendix for a proof of (\ref{rai3})). The
density (\ref{rai}}) \textit{is characterized by fat tails (with polynomial,
instead of exponential, decay) and, as a consequence, the mean waiting time
is infinite.}

\textit{For }$t\rightarrow 0$\textit{\ the density of the interarrival times
displays the following behavior:}%
\begin{equation}
\Pr \left\{ \mathcal{U}_{j}\in dt\right\} /dt\sim \frac{\lambda t^{\nu -1}}{%
\Gamma (\nu )},  \label{rairai}
\end{equation}%
\textit{which means that }$\mathcal{U}_{j}$\textit{\ takes small values with
large probability. Therefore, by considering (\ref{rai}) and (\ref{rairai})
together, we can draw the conclusion that the behavior of the density of the
interarrival times differs from standard Poisson in that the intermediate
values are assumed with smaller probability than in the exponential case.}

\

\noindent \textbf{Remark 2.5}\textit{\ The distribution function (\ref{fk5})
of the waiting-time of the }$k$\textit{-th event coincides, for }$\lambda
=1, $\textit{\ with the so-called Mittag-Leffler distribution of [13]-[19]
(see formula (1) of [13])}%
\begin{equation}
F_{\nu ,t}(x)=\sum_{j=0}^{\infty }\frac{\Gamma (t+j)x^{\nu (t+j)}}{j!\Gamma
(t)\Gamma (1+\nu (t+j))},\qquad x>0,  \label{fk6}
\end{equation}%
\textit{\ when the time argument $t$ takes integer values }$k\geq1.$\textit{%
\ On the contrary the space argument }$x$\textit{\ coincides, in our case,
with the time }$t.$\textit{\ Therefore the process }$X_{\nu }(t)$ \textit{%
with distribution (}\ref{fk6})\textit{\ can be considered as the
continuous-time analogue of the process representing the instant of the }$k$%
\textit{-th event of the fractional Poisson process }$\mathcal{N}_{\nu }.$

\

\noindent \textbf{Remark 2.6} \textit{We observe that also for the
waiting-time density (\ref{fk}) we can find a link with the solution to the
fractional diffusion equation (\ref{due.16}). This can be shown by rewriting
its Laplace transform (\ref{fk2}) as}%
\begin{eqnarray}
\mathcal{L}\left\{ f_{k}^{\nu }(t);s\right\} &=&\frac{\lambda ^{k}}{(s^{\nu
}+\lambda )^{k}}=\mathcal{L}\left\{ f_{k}(t);s^{\nu }\right\}  \notag \\
&=&\int_{0}^{+\infty }e^{-s^{\nu }t}\frac{\lambda ^{k}t^{k-1}}{(k-1)!}%
e^{-\lambda t}dt.  \notag
\end{eqnarray}%
\textit{By using again (\ref{es}) we get}%
\begin{eqnarray}
f_{k}^{\nu }(t) &=&\int_{0}^{+\infty }g_{\nu }(t;y)f_{k}(y)dy  \label{fk4} \\
&=&\int_{0}^{+\infty }g_{\nu }(t;\frac{y}{\lambda })\frac{y^{k-1}e^{-y}}{%
(k-1)!}dy.  \notag
\end{eqnarray}%
\textit{Formula (\ref{fk4}) permits us to conclude that }$f_{k}^{\nu }(t)$%
\textit{\ can be interpreted as a stable law }$S_{\nu }$\textit{\ with a
random scale parameter possessing an Erlang distribution.}

\subsection{The probability generating function}

We consider now the equation governing the probability generating function,
defined, for any $0<u\leq 1$, as%
\begin{equation}
G_{\nu }(u,t)=\sum_{k=0}^{\infty }u^{k}p_{k}^{\nu }(t).  \label{pgf}
\end{equation}%
From (\ref{uno.2}) it is straightforward that it coincides with the solution
to the fractional differential equation%
\begin{equation}
\frac{\partial ^{\nu }G(u,t)}{\partial t^{\nu }}=\lambda (u-1)G(u,t),\qquad
0<\nu \leq 1  \label{due.2}
\end{equation}%
subject to the initial condition $G(u,0)=1.$ As already proved in [2] the
Laplace transform of $G_{\nu }=G_{\nu }(u,t)$ is given by%
\begin{eqnarray}
\mathcal{L}\left\{ G_{\nu }(u,t);s\right\} &=&\int_{0}^{+\infty
}e^{-st}G_{\nu }(u,t)dt  \label{due.3} \\
&=&\frac{s^{\nu -1}}{s^{\nu }-\lambda (u-1)}  \notag
\end{eqnarray}%
so that the probability generating function can be expressed as%
\begin{equation}
G_{\nu }(u,t)=E_{\nu ,1}(\lambda (u-1)t^{\nu }).  \label{due.4}
\end{equation}%
By considering (\ref{due.4}) together with the previous results we get the
following relationship, valid for the infinite sum of GML functions:%
\begin{equation}
\sum_{k=0}^{\infty }(\lambda ut^{\nu })^{k}E_{\nu ,\nu k+1}^{k+1}(-\lambda
t^{\nu })=E_{\nu ,1}(\lambda (u-1)t^{\nu }).  \label{rel}
\end{equation}%
For $u=1$ it shows again that $\sum_{k=0}^{\infty }p_{k}^{\nu }(t)=1.$ The
result (\ref{rel}) can be checked by resorting to the Laplace transforms and
noting that%
\begin{eqnarray}
\sum_{k=0}^{\infty }u^{k}\mathcal{L}\left\{ p_{k}^{\nu }(t);s\right\}
&=&s^{\nu -1}\sum_{k=0}^{\infty }\frac{(u\lambda )^{k}}{(s^{\nu }+\lambda
)^{k+1}}  \label{rel2} \\
&=&\frac{s^{\nu -1}}{s^{\nu }-\lambda (u-1)}.  \notag
\end{eqnarray}%
Formula (\ref{rel}) suggests a useful general relationship between the
infinite sum of GML functions and the standard Mittag-Leffler function:%
\begin{equation}
\sum_{k=0}^{\infty }(ux)^{k}E_{\nu ,\nu k+1}^{k+1}(-x)=E_{\nu
,1}(x(u-1)),\quad 0<u\leq 1.  \label{rel3}
\end{equation}

For $u=1$ it shows again that $\sum_{k=0}^{\infty }p_{k}^{\nu }(t)=1.$

By considering the derivatives of the probability generating function (\ref%
{due.4}) we can easily derive the factorial moments of $\mathbb{\mathcal{N}}%
_{\nu }(t)$ which read%
\begin{equation}
\mathbb{E}\left[ \mathbb{\mathcal{N}}_{\nu }(t)(\mathbb{\mathcal{N}}_{\nu
}(t)-1)...(\mathbb{\mathcal{N}}_{\nu }(t)-r+1)\right] =\frac{\left( \lambda
t^{\nu }\right) ^{r}r!}{\Gamma (\nu r+1)}.  \label{rel6}
\end{equation}%
These are particularly useful in checking that $\mathbb{\mathcal{N}}_{\nu
}(t),t>0$ represents a Cox process with directing measure $\Lambda .$
Indeed, as pointed out in [11], the factorial moments of a Cox process
coincide with the ordinary moments of its directing measure. We show that
this holds for $\mathbb{\mathcal{N}}_{\nu }$, by using the contour integral
representation of the inverse of Gamma function,%
\begin{eqnarray*}
\mathbb{E}\left[ \Lambda (\left( 0,t\right] )\right] ^{r}
&=&\int_{0}^{+\infty }y^{r}v_{2\nu }(y,t)dy \\
&=&\int_{0}^{+\infty }\frac{y^{r}}{\lambda t^{\nu }}W_{-\nu ,1-\nu }\left( -%
\frac{y}{\lambda t^{\nu }}\right) dy \\
&=&\frac{1}{\lambda t^{\nu }}\frac{1}{2\pi i}\int_{0}^{+\infty
}y^{r}dy\int_{Ha}\frac{e^{z-\frac{yt^{-\nu }}{\lambda }z^{\nu }}}{z^{1-\nu }}%
dz \\
&=&\frac{\lambda ^{r}}{2\pi i}\int_{Ha}\frac{e^{z}}{z^{1+\nu r}}%
dz\int_{0}^{+\infty }t^{\nu r}w^{r}e^{-w}dw \\
&=&\frac{\lambda ^{r}t^{\nu r}}{2\pi i}\Gamma (r+1)\int_{Ha}\frac{e^{z}}{%
z^{1+\nu r}}dz=\frac{\lambda ^{r}t^{\nu r}r!}{\Gamma (\nu r+1)},
\end{eqnarray*}%
which coincides with (\ref{rel6}).

\

By applying formula (\ref{rel3}) it is easy to obtain also the moments
generating function of the distribution, defined as%
\begin{equation}
M_{\nu }(t,\mu )=\sum_{k=0}^{\infty }e^{-\mu k}p_{k}^{\nu }(t),\qquad \mu >0,
\label{rel4}
\end{equation}%
which is given by%
\begin{eqnarray}
M_{\nu }(t,\mu ) &=&\sum_{k=0}^{\infty }(\lambda e^{-\mu }t^{\nu
})^{k}E_{\nu ,\nu k+1}^{k+1}(-\lambda t^{\nu })  \label{rel5} \\
&=&E_{\nu ,1}(\lambda (e^{-\mu }-1)t^{\nu }).  \notag
\end{eqnarray}%
The $r$-th moments of the distribution can be obtained by successively
deriving (\ref{rel4}) as%
\begin{eqnarray*}
\mathbb{E\mathcal{N}}_{\nu }(t)^{r} &\mathbb{=}&\mathbb{(-}1)^{r}\left[
\frac{d^{r}}{d\mu ^{r}}M_{\nu }(t,\mu )\right] _{\mu =0} \\
&=&\mathbb{(-}1)^{r}\left[ \frac{d^{r}}{d\mu ^{r}}E_{\nu ,1}(\lambda
(e^{-\mu }-1)t^{\nu })\right] _{\mu =0} \\
&=&\sum_{j=1}^{r}\frac{C_{j,r}\left( \lambda t^{\nu }\right) ^{j}}{\Gamma
(\nu j+1)},
\end{eqnarray*}%
where $C_{j,r}$ are constants. For $r=1$ we get%
\begin{eqnarray}
m_{\nu }(t) &=&\mathbb{E\mathcal{N}}_{\nu }(t)\mathbb{=}\mathbb{-}\lambda
t^{\nu }\left[ \frac{d}{d\mu }E_{\nu ,1}(\lambda (e^{-\mu }-1)t^{\nu })%
\right] _{\mu =0}  \label{ren} \\
&=&\lambda t^{\nu }\left[ \frac{d}{dx}E_{\nu ,1}(x)\right] _{x=0}=\frac{%
\lambda t^{\nu }}{\Gamma (\nu +1)},  \notag
\end{eqnarray}%
which coincides with the renewal function evaluated in [15], for $\lambda =1$%
. It is evident also from (\ref{ren}) that the mean waiting time (which
coincides with $\lim_{t\rightarrow \infty }t/m_{\nu }(t)$) is infinite,
since $\nu <1.$

\section{Second-type fractional recursive differential equation}

\subsection{The solution}

In this section we generalize the results obtained so far to a fractional
recursive differential equation containing two time-fractional derivatives.
We show that some properties of the first model of fractional Poisson
process are still valid: the solutions represent, for $k\geq 0,$ a proper
probability distribution and the corresponding process is again a renewal
process. Moreover the density of the interarrival times display the same
asymptotic behavior of the previous model.

We consider the following recursive differential equation%
\begin{equation}
\frac{d^{2\nu }p_{k}}{dt^{2\nu }}+2\lambda \frac{d^{\nu }p_{k}}{dt^{\nu }}%
=-\lambda ^{2}(p_{k}-p_{k-1}),\quad k\geq 0,  \label{sec1}
\end{equation}%
where $\nu \in \left( 0,1\right] ,$ subject to the initial conditions%
\begin{eqnarray}
p_{k}(0) &=&\left\{
\begin{array}{c}
1\qquad k=0 \\
0\qquad k\geq 1%
\end{array}%
\right. ,\quad \text{for }0<\nu \leq 1  \label{sec2} \\
p_{k}^{\prime }(0) &=&0,\qquad k\geq 0,\quad \text{for }\frac{1}{2}<\nu \leq
1  \notag
\end{eqnarray}%
and $p_{-1}(t)=0$. In the following theorem we derive the solution to (\ref%
{sec1})-(\ref{sec2}), which can be still expressed in terms of GML functions.

\

\noindent \textbf{Theorem 3.1 }\textit{The solution }$\widehat{p}_{k}^{\nu
}(t)$\textit{, for }$k=0,1,...$\textit{\ and }$t\geq 0,$\textit{\ of the
Cauchy problem (\ref{sec1})-(\ref{sec2}) is given by}%
\begin{equation}
\widehat{p}_{k}^{\nu }(t)=\lambda ^{2k}t^{2k\nu }E_{\nu ,2k\nu
+1}^{2k+1}(-\lambda t^{\nu })+\lambda ^{2k+1}t^{(2k+1)\nu }E_{\nu ,(2k+1)\nu
+1}^{2k+2}(-\lambda t^{\nu }),\quad \text{ }k\geq 0,t>0.  \label{sec3}
\end{equation}%
\textbf{Proof }Following the lines of the proof of Theorem 2.1, we take the
Laplace transform of equation (\ref{sec1}) together with the conditions (\ref%
{sec2}), thus obtaining the following recursive formula, for $k\geq 1$%
\begin{eqnarray}
\mathcal{L}\left\{ \widehat{p}_{k}^{\nu }(t);s\right\} &=&\frac{\lambda ^{2}%
}{s^{2\nu }+2\lambda s^{\nu }+\lambda ^{2}}\mathcal{L}\left\{ \widehat{p}%
_{k-1}^{\nu }(t);s\right\} \\
&=&\frac{\lambda ^{2}}{(s^{\nu }+\lambda )^{2}}\mathcal{L}\left\{ \widehat{p}%
_{k-1}^{\nu }(t);s\right\} ,  \notag
\end{eqnarray}%
while, for $k=0$, we get
\begin{equation}
\mathcal{L}\left\{ \widehat{p}_{0}^{\nu }(t);s\right\} =\frac{s^{2\nu
-1}+2\lambda s^{\nu -1}}{s^{2\nu }+2\lambda s^{\nu }+\lambda ^{2}}.
\end{equation}%
Therefore the Laplace transform of the solution reads%
\begin{equation}
\mathcal{L}\left\{ \widehat{p}_{k}^{\nu }(t);s\right\} =\frac{\lambda
^{2k}s^{2\nu -1}+2\lambda ^{2k+1}s^{\nu -1}}{\left( s^{\nu }+\lambda \right)
^{2k+2}}.  \label{lap}
\end{equation}%
We can invert (\ref{lap}) by using (\ref{pra}) with $\delta =2k+2,$ $\beta
=\nu $ and $\gamma =2k\nu +1$ or $\gamma =(2k+1)\nu +1,$ thus obtaining the
following expression%
\begin{equation}
\widehat{p}_{k}^{\nu }(t)=\lambda ^{2k}t^{2\nu k}E_{\nu ,2k\nu
+1}^{2k+2}(-\lambda t^{\nu })+2\lambda ^{2k+1}t^{(2k+1)\nu }E_{\nu
,(2k+1)\nu +1}^{2k+2}(-\lambda t^{\nu }).  \label{pre}
\end{equation}%
We prove now the following general formula holding for a sum of GML
functions:

\begin{equation}
x^{n}E_{\nu ,n\nu +z}^{m}(-x)+x^{n+1}E_{\nu ,(n+1)\nu
+z}^{m}(-x)=x^{n}E_{\nu ,n\nu +z}^{m-1}(-x),\quad \text{ }n,m>0,z\geq 0,x>0,
\label{gen}
\end{equation}%
which can be proved by rewriting the l.h.s. as follows:%
\begin{eqnarray*}
&&\frac{x^{n}}{(m-1)!}\sum_{j=0}^{\infty }\frac{(m-1+j)!(-x)^{j}}{j!\Gamma
(\nu j+n\nu +z)}-\frac{x^{n}}{(m-1)!}\sum_{j=0}^{\infty }\frac{%
(m-1+j)!(-x)^{j+1}}{j!\Gamma (\nu j+(n+1)\nu +z)} \\
&=&\frac{x^{n}}{(m-1)!}\sum_{j=0}^{\infty }\frac{(m-1+j)!(-x)^{j}}{j!\Gamma
(\nu j+n\nu +z)}-\frac{x^{n}}{(m-1)!}\sum_{l=1}^{\infty }\frac{%
(m+l-2)!(-x)^{l}}{(l-1)!\Gamma (\nu l+n\nu +z)} \\
&=&\frac{x^{n}}{(m-1)!}\sum_{l=1}^{\infty }\frac{(m+l-2)!(-x)^{l}}{%
(l-1)!\Gamma (\nu l+n\nu +z)}\left[ \frac{m-1+l}{l}-1\right] +\frac{x^{n}}{%
\Gamma (n\nu +z)} \\
&=&\frac{x^{n}}{(m-2)!}\sum_{l=1}^{\infty }\frac{(m+l-2)!(-x)^{l}}{l!\Gamma
(\nu l+n\nu +z)}+\frac{x^{n}}{\Gamma (n\nu +z)}=x^{n}E_{\nu ,n\nu
+z}^{m-1}(-x)
\end{eqnarray*}%
For $m=2k+2,$ $z=1,x=\lambda t^{\nu }$ and $n=2k$ formula (\ref{gen}) gives
the following identity:%
\begin{equation*}
\lambda ^{2k}t^{2\nu k}E_{\nu ,2k\nu +1}^{2k+2}(-\lambda t^{\nu })+\lambda
^{2k+1}t^{(2k+1)\nu }E_{\nu ,(2k+1)\nu +1}^{2k+2}(-\lambda t^{\nu })=\lambda
^{2k}t^{2\nu k}E_{\nu ,2k\nu +1}^{2k+1}(-\lambda t^{\nu }),
\end{equation*}%
which coincides with the first term in (\ref{sec3}).

It remains to check only that the initial conditions in (\ref{sec2}) hold:
the first one is clearly satisfied since it is, for $k=0,$%
\begin{eqnarray*}
\widehat{p}_{0}^{\nu }(t) &=&\sum_{r=0}^{\infty }\frac{(-\lambda )^{r}t^{\nu
r}}{\Gamma (\nu r+1)}+\lambda \sum_{r=0}^{\infty }\frac{(r+1)(-\lambda
)^{r}t^{\nu (r+1)}}{\Gamma (\nu r+\nu +1)} \\
&=&\left[ \text{for }t=0\right] =1
\end{eqnarray*}%
and, for $k\geq 1$,
\begin{equation*}
\widehat{p}_{k}^{\nu }(t)=\frac{\lambda ^{2k}}{(2k)!}\sum_{r=0}^{\infty }%
\frac{(2k+r)!(-\lambda )^{r}t^{\nu (2k+r)}}{r!\Gamma (\nu r+2k\nu +1)}+\frac{%
\lambda ^{2k+1}}{(2k+1)!}\sum_{r=0}^{\infty }\frac{(2k+r+1)!(-\lambda
)^{r}t^{\nu (2k+r+1)}}{r!\Gamma (\nu r+2k\nu +\nu +1)},
\end{equation*}%
which vanishes for $t=0.$ The second condition in (\ref{sec2}) is
immediately verified for $k\geq 1$, since it is%
\begin{equation}
\frac{d}{dt}\widehat{p}_{k}^{\nu }(t)=\frac{\lambda ^{2k}}{(2k)!}%
\sum_{r=1}^{\infty }\frac{(2k+r)!(-\lambda )^{r}t^{\nu (2k+r)-1}}{r!\Gamma
(\nu r+2k\nu )}+\frac{\lambda ^{2k+1}}{(2k+1)!}\sum_{r=0}^{\infty }\frac{%
(2k+r+1)!(-\lambda )^{r}t^{\nu (2k+r+1)-1}}{r!\Gamma (\nu r+2k\nu +\nu )},
\label{che}
\end{equation}%
which for $t=0$ vanishes in the interval $\frac{1}{2}<\nu \leq 1.$ Then we
check that this happens also for $k=0$: indeed in this case (\ref{che})
reduces to
\begin{eqnarray*}
\frac{d}{dt}\widehat{p}_{0}^{\nu }(t) &=&\sum_{r=1}^{\infty }\frac{(-\lambda
)^{r}t^{\nu r-1}}{\Gamma (\nu r)}+\lambda \sum_{r=0}^{\infty }\frac{%
(r+1)^{2}(-\lambda )^{r}t^{\nu (r+1)-1}}{\Gamma (\nu r+\nu )} \\
&=&\sum_{r=2}^{\infty }\frac{(-\lambda )^{r}t^{\nu r-1}}{\Gamma (\nu r)}-%
\frac{\lambda t^{\nu -1}}{\Gamma (\nu )}+\lambda \sum_{r=1}^{\infty }\frac{%
(r+1)^{2}(-\lambda )^{r}t^{\nu (r+1)-1}}{\Gamma (\nu r+\nu )}+\frac{\lambda
t^{\nu -1}}{\Gamma (\nu )} \\
&=&\left[ \text{for }t=0\right] =0.
\end{eqnarray*}%
\hfill $\square $

\

\noindent \textbf{Remark 3.1}\textit{\ The solution (\ref{sec3}) can be
expressed in terms of the solution (\ref{gml3}) of the first model as follows%
}%
\begin{equation}
\widehat{p}_{k}^{\nu }(t)=p_{2k}^{\nu }(t)+p_{2k+1}^{\nu }(t).  \label{sec4}
\end{equation}%
\textit{Therefore it can be interpreted, for }$k=0,1,2,...,$\textit{\ as}
\textit{the probability distribution }$\Pr \left\{ \widehat{\mathcal{N}}%
_{\nu }(t)=k\right\} $\textit{\ for a process }$\widehat{\mathcal{N}}_{\nu
}. $ \textit{Indeed, by (\ref{sec4}), we get}%
\begin{equation}
\Pr \left\{ \widehat{\mathcal{N}}_{\nu }(t)=k\right\} =\Pr \left\{ \mathcal{N%
}_{\nu }(t)=2k\right\} +\Pr \left\{ \mathcal{N}_{\nu }(t)=2k+1\right\} ,
\label{sec9}
\end{equation}%
\textit{so that it is}%
\begin{eqnarray}
&&\sum_{k=0}^{\infty }\Pr \left\{ \widehat{\mathcal{N}}_{\nu }(t)=k\right\}
\label{ch} \\
&=&\sum_{k=0}^{\infty }\Pr \left\{ \mathcal{N}_{\nu }(t)=2k\right\}
+\sum_{k=0}^{\infty }\Pr \left\{ \mathcal{N}_{\nu }(t)=2k+1\right\} =1
\notag
\end{eqnarray}%
\textit{Moreover the relationship (\ref{sec4}) shows that the process
governed by the second-type equation can be seen as a first-type fractional
process, which jumps upward at even-order events }$A_{2k}$\textit{\ while
the probability of the successive odd-indexed events }$A_{2k+1}$\textit{\ is
added to that of }$A_{2k}$\textit{.}

\

We can check that expression (\ref{sec4}) is the solution to equation (\ref%
{sec1}), subject to the initial conditions (\ref{sec2}), by using the form
of $p_{k}^{\nu }$\ appearing in the last line of (\ref{rec3}) which is more
suitable to this aim:%
\begin{equation}
\widehat{p}_{k}^{\nu }(t)=\sum_{r=2k}^{\infty }\binom{r}{2k}\frac{(-\lambda
t^{\nu })^{r}}{\Gamma \left( \nu r+1\right) }-\sum_{r=2k+1}^{\infty }\binom{r%
}{2k+1}\frac{(-\lambda t^{\nu })^{r}}{\Gamma \left( \nu r+1\right) }.
\label{der}
\end{equation}

The fractional derivatives of (\ref{der}) can be evaluated by applying the
definition (\ref{dc}), as follows:%
\begin{eqnarray}
\frac{d^{\nu }}{dt^{\nu }}\widehat{p}_{k}^{\nu }(t) &=&\frac{1}{\Gamma
(1-\nu )}\sum_{r=2k}^{\infty }\binom{r}{2k}\frac{(-\lambda )^{r}\nu r}{%
\Gamma \left( \nu r+1\right) }\int_{0}^{t}\frac{s^{\nu r-1}}{(t-s)^{\nu }}ds+
\label{der1} \\
&&-\frac{1}{\Gamma (1-\nu )}\sum_{r=2k+1}^{\infty }\binom{r}{2k+1}\frac{%
(-\lambda )^{r}\nu r}{\Gamma \left( \nu r+1\right) }\int_{0}^{t}\frac{s^{\nu
r-1}}{(t-s)^{\nu }}ds  \notag \\
&=&\frac{1}{\Gamma (1-\nu )}\sum_{r=2k}^{\infty }\binom{r}{2k}\frac{\Gamma
(\nu r+1)\Gamma (1-\nu )}{\Gamma (\nu r-\nu +1)\Gamma \left( \nu r+1\right) }%
(-\lambda )^{r}t^{\nu r-\nu }+  \notag \\
&&-\frac{1}{\Gamma (1-\nu )}\sum_{r=2k+1}^{\infty }\binom{r}{2k+1}\frac{%
\Gamma (\nu r+1)\Gamma (1-\nu )}{\Gamma (\nu r-\nu +1)\Gamma \left( \nu
r+1\right) }(-\lambda )^{r}t^{\nu r-\nu }  \notag
\end{eqnarray}%
and, for $\nu \in \left[ 0,\frac{1}{2}\right] $\ and $m=1,$%
\begin{eqnarray}
\frac{d^{2\nu }}{dt^{2\nu }}\widehat{p}_{k}^{\nu }(t) &=&\frac{1}{\Gamma
(1-2\nu )}\sum_{r=2k}^{\infty }\binom{r}{2k}\frac{(-\lambda )^{r}\nu r}{%
\Gamma \left( \nu r+1\right) }\int_{0}^{t}\frac{s^{\nu r-1}}{(t-s)^{2\nu }}%
ds+  \label{der2} \\
&&-\frac{1}{\Gamma (1-2\nu )}\sum_{r=2k+1}^{\infty }\binom{r}{2k+1}\frac{%
(-\lambda )^{r}\nu r}{\Gamma \left( \nu r+1\right) }\int_{0}^{t}\frac{s^{\nu
r-1}}{(t-s)^{2\nu }}ds  \notag \\
&=&\frac{1}{\Gamma (1-2\nu )}\sum_{r=2k}^{\infty }\binom{r}{2k}\frac{\Gamma
(\nu r+1)\Gamma (1-2\nu )}{\Gamma (\nu r-2\nu +1)\Gamma \left( \nu
r+1\right) }(-\lambda )^{r}t^{\nu r-2\nu }+  \notag \\
&&-\frac{1}{\Gamma (1-2\nu )}\sum_{r=2k+1}^{\infty }\binom{r}{2k+1}\frac{%
\Gamma (\nu r+1)\Gamma (1-2\nu )}{\Gamma (\nu r-\nu +1)\Gamma \left( \nu
r+1\right) }(-\lambda )^{r}t^{\nu r-2\nu }.  \notag
\end{eqnarray}%
It is easy to check that the last expression is valid also in the case $\nu
\in \left[ \frac{1}{2},1\right] $\ and $m=2.$\ By considering (\ref{der1})
and (\ref{der2}) together we get%
\begin{eqnarray}
&&\frac{d^{2\nu }}{dt^{2\nu }}\widehat{p}_{k}^{\nu }+2\lambda \frac{d^{\nu }%
}{dt^{\nu }}\widehat{p}_{k}^{\nu }  \label{der3} \\
&=&\sum_{r=2k}^{\infty }\binom{r}{2k}\frac{(-\lambda )^{r}t^{\nu r-2\nu }}{%
\Gamma (\nu r-2\nu +1)}-\sum_{r=2k+1}^{\infty }\binom{r}{2k+1}\frac{%
(-\lambda )^{r}t^{\nu r-2\nu }}{\Gamma (\nu r-2\nu +1)}+  \notag \\
&&+2\lambda \left[ \sum_{r=2k}^{\infty }\binom{r}{2k}\frac{(-\lambda
)^{r}t^{\nu r-\nu }}{\Gamma (\nu r-\nu +1)}-\sum_{r=2k+1}^{\infty }\binom{r}{%
2k+1}\frac{(-\lambda )^{r}t^{\nu r-\nu }}{\Gamma (\nu r-\nu +1)}\right]
\notag \\
&=&\sum_{r=2k-2}^{\infty }\binom{r+2}{2k}\frac{(-1)^{r}\lambda ^{r+2}t^{\nu
r}}{\Gamma (\nu r+1)}-\sum_{r=2k-1}^{\infty }\binom{r+2}{2k+1}\frac{%
(-1)^{r}\lambda ^{r+2}t^{\nu r}}{\Gamma (\nu r+1)}+  \notag \\
&&+2\lambda \left[ -\sum_{r=2k-1}^{\infty }\binom{r+1}{2k}\frac{%
(-1)^{r}\lambda ^{r+1}t^{\nu r}}{\Gamma (\nu r+1)}+\sum_{r=2k}^{\infty }%
\binom{r+1}{2k+1}\frac{(-1)^{r}\lambda ^{r+1}t^{\nu r}}{\Gamma (\nu r+1)}%
\right]  \notag \\
&=&\sum_{r=2k-2}^{\infty }\binom{r+2}{2k}A_{r}-\sum_{r=2k-1}^{\infty }\binom{%
r+2}{2k+1}A_{r}+  \notag \\
&&-2\sum_{r=2k-1}^{\infty }\binom{r+1}{2k}A_{r}+2\sum_{r=2k}^{\infty }\binom{%
r+1}{2k+1}A_{r}  \notag
\end{eqnarray}%
for $A_{r}=\frac{(-1)^{r}\lambda ^{r+2}t^{\nu r}}{\Gamma (\nu r+1)}$\ where
the second step follows by putting $r^{\prime }=r-2$\ in the first two sums
and $r^{\prime }=r-1$\ in the second ones.

We want to show that (\ref{der3}) is equal to%
\begin{eqnarray}
&&-\lambda ^{2}(\widehat{p}_{k}^{\nu }-\widehat{p}_{k-1}^{\nu })
\label{der4} \\
&=&-\lambda ^{2}(p_{2k}^{\nu }+p_{2k+1}^{\nu }-p_{2k-2}^{\nu }-p_{2k-1}^{\nu
})  \notag \\
&=&-\lambda ^{2}\left[ \sum_{r=2k}^{\infty }\binom{r}{2k}\frac{(-\lambda
)^{r}t^{\nu r}}{\Gamma (\nu r+1)}-\sum_{r=2k+1}^{\infty }\binom{r}{2k+1}%
\frac{(-\lambda )^{r}t^{\nu r}}{\Gamma (\nu r+1)}\right. +  \notag \\
&&\left. -\sum_{r=2k-2}^{\infty }\binom{r}{2k-2}\frac{(-\lambda )^{r}t^{\nu
r}}{\Gamma (\nu r+1)}+\sum_{r=2k-1}^{\infty }\binom{r}{2k-1}\frac{(-\lambda
)^{r}t^{\nu r}}{\Gamma (\nu r+1)}\right]  \notag \\
&=&-\sum_{r=2k}^{\infty }\binom{r}{2k}A_{r}+\sum_{r=2k+1}^{\infty }\binom{r}{%
2k+1}A_{r}+  \notag \\
&&+\sum_{r=2k-2}^{\infty }\binom{r}{2k-2}A_{r}-\sum_{r=2k-1}^{\infty }\binom{%
r}{2k-1}A_{r}.  \notag
\end{eqnarray}%
Then, by considering together (\ref{der3}) and (\ref{der4}), we get%
\begin{eqnarray}
&&\frac{d^{2\nu }}{dt^{2\nu }}\widehat{p}_{k}^{\nu }+2\lambda \frac{d^{\nu }%
}{dt^{\nu }}\widehat{p}_{k}^{\nu }+\lambda ^{2}(\widehat{p}_{k}^{\nu }-%
\widehat{p}_{k-1}^{\nu })  \label{der5} \\
&=&\sum_{r=2k-2}^{\infty }\left[ \binom{r+2}{2k}-\binom{r}{2k-2}\right]
A_{r}-\sum_{r=2k-1}^{\infty }\left[ \binom{r+2}{2k+1}+2\binom{r+1}{2k}-%
\binom{r}{2k-1}\right] A_{r}+  \notag \\
&&+\sum_{r=2k}^{\infty }\left[ 2\binom{r+1}{2k+1}+\binom{r}{2k}\right]
A_{r}-\sum_{r=2k+1}^{\infty }\binom{r}{2k+1}A_{r}  \notag \\
&=&\sum_{r=2k+1}^{\infty }\left[ \binom{r+2}{2k}-2\binom{r+1}{2k}+\binom{r}{%
2k-1}-\binom{r+2}{2k+1}+2\binom{r+1}{2k+1}\right. +  \notag \\
&&\left. -\binom{r}{2k-2}+\binom{r}{2k}-\binom{r}{2k+1}\right] A_{r}+\mathbf{%
R}  \notag
\end{eqnarray}%
where%
\begin{eqnarray}
\mathbf{R} &=&A_{2k-2}\left[ 1-1\right] +A_{2k-1}\left[ \binom{2k+1}{2k}-%
\binom{2k-1}{2k-2}\right] +  \label{der6} \\
&&+A_{2k}\left[ \binom{2k+2}{2k}-\binom{2k}{2k-2}\right] -A_{2k-1}\left[
\binom{2k+1}{2k+1}+2\binom{2k}{2k}-\binom{2k-1}{2k-1}\right] +  \notag \\
&&-A_{2k}\left[ \binom{2k+2}{2k+1}+2\binom{2k+1}{2k}-\binom{2k}{2k-1}\right]
+A_{2k}\left[ 2\binom{2k+1}{2k+1}+\binom{2k}{2k}\right]  \notag \\
&=&0.  \notag
\end{eqnarray}%
The sum appearing in (\ref{der5}) can be developed by considering that%
\begin{eqnarray*}
\binom{r+2}{2k}-\binom{r+1}{2k} &=&\frac{(r+1)!}{(2k)!(r+1-2k)!}\left[ \frac{%
r+2}{r+2-2k}-1\right] \\
&=&\binom{r+1}{2k-1}
\end{eqnarray*}%
and analogously%
\begin{equation*}
\binom{r}{2k-1}-\binom{r+1}{2k}=-\binom{r}{2k}
\end{equation*}%
\begin{equation*}
-\binom{r+2}{2k+1}+\binom{r+1}{2k+1}=-\binom{r+1}{2k},
\end{equation*}%
so that, by considering (\ref{der6}), we can rewrite (\ref{der5}) as follows%
\begin{eqnarray}
&&\sum_{r=2k+1}^{\infty }\left[ \binom{r+1}{2k-1}-\binom{r}{2k-2}-\binom{r+1%
}{2k}+\binom{r+1}{2k+1}-\binom{r}{2k}\right. +  \label{der8} \\
&&\left. +\binom{r}{2k}-\binom{r}{2k+1}\right] A_{r}  \notag \\
&=&\sum_{r=2k+1}^{\infty }\left[ \binom{r}{2k-1}-\binom{r+1}{2k}+\binom{r+1}{%
2k+1}-\binom{r}{2k+1}\right] A_{r},  \notag
\end{eqnarray}%
since%
\begin{equation}
\binom{r+1}{2k-1}-\binom{r}{2k-2}=\binom{r}{2k-1}.  \label{der7}
\end{equation}%
If we consider now the first two terms of (\ref{der8}) we get%
\begin{equation*}
\binom{r}{2k-1}-\binom{r+1}{2k}=-\binom{r}{2k}
\end{equation*}%
so that we get%
\begin{equation*}
\sum_{r=2k+1}^{\infty }\left[ -\binom{r}{2k}+\binom{r+1}{2k+1}-\binom{r}{2k+1%
}\right] A_{r}=0,
\end{equation*}%
since%
\begin{equation*}
\binom{r+1}{2k+1}-\binom{r}{2k+1}=\binom{r}{2k}.
\end{equation*}

\subsection{The probability generating function}

As we did for the first model we evaluate the probability generating
function and we show that it coincides with the solution to a fractional
equation which arises in the study of the fractional telegraph process (see
[17]).

\

\noindent \textbf{Theorem 3.2 }The probability generating function $\widehat{%
G}_{\nu }(u,t)=\sum_{k=0}^{\infty }u^{k}\widehat{p}_{k}^{\nu }(t),$ $0<u\leq
1,$ coincides with the solution to the following fractional differential
equation%
\begin{equation}
\frac{\partial ^{2\nu }G(u,t)}{\partial t^{2\nu }}+2\lambda \frac{\partial
^{\nu }G(u,t)}{\partial t^{\nu }}=\lambda ^{2}(u-1)G(u,t),\qquad 0<\nu \leq 1
\label{sec5}
\end{equation}%
subject to the initial condition $G(u,0)=1$ and the additional condition $%
G_{t}(u,0)=0$ for $1/2<\nu <1.$ The explicit expression is given by%
\begin{equation}
\widehat{G}_{\nu }(u,t)=\frac{\sqrt{u}+1}{2\sqrt{u}}E_{\nu ,1}(-\lambda (1-%
\sqrt{u})t^{\nu })+\frac{\sqrt{u}-1}{2\sqrt{u}}E_{\nu ,1}(-\lambda (1+\sqrt{u%
})t^{\nu }).  \label{sec6}
\end{equation}

\noindent \textbf{Proof }By applying the Laplace transform to (\ref{sec5}),
we get%
\begin{equation*}
(s^{2\nu }+2\lambda s^{\nu })\mathcal{L(}\widehat{G}_{\nu }(u,t);s)+(s^{2\nu
-1}+2\lambda s^{\nu -1})=\lambda ^{2}(u-1)\mathcal{L(}\widehat{G}_{\nu
}(u,t);s)
\end{equation*}%
and then%
\begin{equation}
\mathcal{L(}\widehat{G}_{\nu }(u,t);s)=\frac{s^{2\nu -1}+2\lambda s^{\nu -1}%
}{s^{2\nu }+2\lambda s^{\nu }+\lambda ^{2}(1-u)}.  \label{sec7}
\end{equation}%
We can recognize in (\ref{sec5}) the fractional equation satisfied by the
characteristic function of the fractional telegraph process\ studied in [17]
(see formula (2.3a) with $c^{2}\beta ^{2}=\lambda ^{2}(1-u)$) and thus the
Laplace transform (\ref{sec7}) coincides with formula (2.6) therein. By
applying the result of Theorem 2.1 of the cited paper, we obtain the inverse
Laplace transform of (\ref{sec7}) as given in (\ref{sec6}).\hfill $\square $

\

\noindent \textbf{Remark 3.2}\textit{\ As a first check we note that (\ref%
{sec6}) reduces to one for }$u=1,$ \textit{so that we prove again that (\ref%
{ch}) holds. Moreover, as an alternative proof of the previous theorem, we
can show that the series expansion of (\ref{sec6}) coincides with }$%
\sum_{k=0}^{\infty }u^{k}\widehat{p}_{k}^{\nu }(t)$\textit{\ for }$\widehat{p%
}_{k}^{\nu }(t)$\textit{\ given in (\ref{sec3}):}%
\begin{eqnarray*}
&&\widehat{G}_{\nu }(u,t) \\
&=&\frac{\sqrt{u}+1}{2\sqrt{u}}\sum_{j=0}^{\infty }\frac{(-\lambda (1-\sqrt{u%
})t^{\nu })^{j}}{\Gamma (\nu j+1)}+\frac{\sqrt{u}-1}{2\sqrt{u}}%
\sum_{j=0}^{\infty }\frac{(-\lambda (1+\sqrt{u})t^{\nu })^{j}}{\Gamma (\nu
j+1)} \\
&=&\frac{\sqrt{u}+1}{2\sqrt{u}}\sum_{j=0}^{\infty }\frac{(\lambda t^{\nu
})^{j}}{\Gamma (\nu j+1)}\sum_{k=0}^{j}\binom{j}{k}\left( \sqrt{u}\right)
^{k}(-1)^{j-k}+\frac{\sqrt{u}-1}{2\sqrt{u}}\sum_{j=0}^{\infty }\frac{%
(-\lambda t^{\nu })^{j}}{\Gamma (\nu j+1)}\sum_{k=0}^{j}\binom{j}{k}\left(
\sqrt{u}\right) ^{k} \\
&=&\frac{\sqrt{u}+1}{2\sqrt{u}}\sum_{k=0}^{\infty }\left( -\sqrt{u}\right)
^{k}\sum_{j=k}^{\infty }\frac{(-\lambda t^{\nu })^{j}}{\Gamma (\nu j+1)}%
\binom{j}{k}+\frac{\sqrt{u}-1}{2\sqrt{u}}\sum_{k=0}^{\infty }\left( \sqrt{u}%
\right) ^{k}\sum_{j=k}^{\infty }\frac{(-\lambda t^{\nu })^{j}}{\Gamma (\nu
j+1)}\binom{j}{k}
\end{eqnarray*}%
\begin{eqnarray*}
&=&\frac{\sqrt{u}+1}{2\sqrt{u}}\sum_{k=0}^{\infty }\left( \sqrt{u}\lambda
t^{\nu }\right) ^{k}\sum_{l=0}^{\infty }\frac{(-\lambda t^{\nu })^{l}}{%
\Gamma (\nu l+\nu k+1)}\binom{l+k}{k}+ \\
&&+\frac{\sqrt{u}-1}{2\sqrt{u}}\sum_{k=0}^{\infty }\left( -\sqrt{u}\lambda
t^{\nu }\right) ^{k}\sum_{l=0}^{\infty }\frac{(-\lambda t^{\nu })^{l}}{%
\Gamma (\nu l+\nu k+1)}\binom{l+k}{k} \\
&=&\frac{\sqrt{u}+1}{2\sqrt{u}}\sum_{k=0}^{\infty }\left( \sqrt{u}\lambda
t^{\nu }\right) ^{k}E_{\nu ,\nu k+1}^{k+1}(-\lambda t^{\nu })+\frac{\sqrt{u}%
-1}{2\sqrt{u}}\sum_{k=0}^{\infty }\left( -\sqrt{u}\lambda t^{\nu }\right)
^{k}E_{\nu ,\nu k+1}^{k+1}(-\lambda t^{\nu }) \\
&=&\frac{\sqrt{u}+1}{2\sqrt{u}}\sum_{m=0}^{\infty }\left( \sqrt{u}\lambda
t^{\nu }\right) ^{2m}E_{\nu ,2m\nu +1}^{2m+1}(-\lambda t^{\nu })+\frac{\sqrt{%
u}+1}{2\sqrt{u}}\sum_{m=0}^{\infty }\left( \sqrt{u}\lambda t^{\nu }\right)
^{2m+1}E_{\nu ,(2m+1)\nu +1}^{2m+2}(-\lambda t^{\nu })+ \\
&&+\frac{\sqrt{u}-1}{2\sqrt{u}}\sum_{m=0}^{\infty }\left( \sqrt{u}\lambda
t^{\nu }\right) ^{2m}E_{\nu ,2m\nu +1}^{2m+1}(-\lambda t^{\nu })+\frac{\sqrt{%
u}-1}{2\sqrt{u}}\sum_{m=0}^{\infty }\left( -\sqrt{u}\lambda t^{\nu }\right)
^{2m+1}E_{\nu ,(2m+1)\nu +1}^{2m+2}(-\lambda t^{\nu }) \\
&=&\sum_{m=0}^{\infty }u^{m}\left[ \left( \lambda t^{\nu }\right)
^{2m}E_{\nu ,2m\nu +1}^{2m+1}(-\lambda t^{\nu })+\left( \lambda t^{\nu
}\right) ^{2m+1}E_{\nu ,(2m+1)\nu +1}^{2m+2}(-\lambda t^{\nu })\right] .
\end{eqnarray*}

\subsection{Properties of the corresponding process}

We can prove that $\widehat{\mathcal{N}}_{\nu }(t),t>0$\ represents a
renewal process, by showing that, also for this model, the required
relationship between $\widehat{p}_{k}^{\nu }(t)$\ and distribution function
of the waiting time $\widehat{T}_{k}$\ of the $k$-th event holds:%
\begin{equation}
\widehat{p}_{k}^{\nu }(t)=\Pr \left\{ \widehat{T}_{k}<t\right\} -\Pr \left\{
\widehat{T}_{k+1}<t\right\}  \label{fin}
\end{equation}%
where%
\begin{equation*}
\widehat{T}_{k}=\inf \left\{ t>0:\widehat{\mathcal{N}}_{\nu }(t)=k\right\} .
\end{equation*}%
or alternatively for the Laplace transform of (\ref{fin})%
\begin{equation}
\mathcal{L}\left\{ \widehat{p}_{k}^{\nu }(t);s\right\} =\frac{1}{s}\mathcal{L%
}\left\{ \widehat{f}_{k}^{\nu }(t);s\right\} -\frac{1}{s}\mathcal{L}\left\{
\widehat{f}_{k+1}^{\nu }(t);s\right\} .  \label{lap5}
\end{equation}%
In view of relationship (\ref{sec4}), we can infer that each interarrival
time $\widehat{\mathcal{U}}_{j}$ is distributed as the sum of two
independent interarrival times $\mathcal{U}_{j}$ of the first model and
therefore from (\ref{mailap}) we have that%
\begin{eqnarray}
\int_{0}^{\infty }e^{-st}\Pr \left\{ \widehat{\mathcal{U}}_{j}\in dt\right\}
&=&\left[ \int_{0}^{\infty }e^{-st}\Pr \left\{ \widehat{\mathcal{U}}_{j}\in
dt\right\} \right] ^{2}=\frac{\lambda ^{2}}{\left( s^{\nu }+\lambda \right)
^{2}}  \label{fin2} \\
&=&\mathcal{L}\left\{ \widehat{f}_{1}^{\nu }(t);s\right\} .  \notag
\end{eqnarray}%
By recursively applying (\ref{lap5}), starting with (\ref{fin2}), we arrive
at%
\begin{equation}
\mathcal{L}\left\{ \widehat{f}_{k}^{\nu }(t);s\right\} =\frac{\lambda ^{2k}}{%
(s^{\nu }+\lambda )^{2k}}.  \label{lap7}
\end{equation}%
By applying again (\ref{pra}) for $\beta =\nu ,$ $\gamma =2\nu k$ and $%
\omega =-\lambda $ we invert (\ref{lap7}) and obtain

\begin{equation}
\widehat{f}_{k}^{\nu }(t)=\lambda ^{2k}t^{2\nu k-1}E_{\nu ,2\nu
k}^{2k}(-\lambda t^{\nu }).  \label{lap6}
\end{equation}%
Moreover from (\ref{fin2}) or (\ref{lap6}) we easily get%
\begin{equation}
\Pr \left\{ \widehat{\mathcal{U}}_{j}\in dt\right\} /dt=\mathcal{L}%
^{-1}\left\{ \frac{\lambda ^{2}}{\left( s^{\nu }+\lambda \right) ^{2}}%
;t\right\} =\lambda ^{2}t^{2\nu -1}E_{\nu ,2\nu }^{2}(-\lambda t^{\nu })=%
\widehat{f}_{1}^{\nu }(t).  \label{lap8}
\end{equation}%
Therefore, in this case, both the waiting-time of the $k$-th event and the
interarrival times possess distributions which are expressed in terms of GML
functions.

\

\noindent \textbf{Remark 3.3 }\textit{As a check, we derive directly the
probability density of the }$\widehat{\mathcal{U}}_{j}$\textit{\ from that
of }$\mathcal{U}_{j}$ \textit{given in (\ref{maimai}), without the use of
Laplace transform} \textit{as follows}:%
\begin{eqnarray*}
\int_{0}^{t}f_{1}^{\nu }(s)f_{1}^{\nu }(t-s)ds &=&\lambda
^{2}\int_{0}^{t}s^{\nu -1}E_{\nu ,\nu }(-\lambda s^{\nu })(t-s)^{\nu
-1}E_{\nu ,\nu }(-\lambda (t-s)^{\nu })ds \\
&=&\lambda ^{2}\sum_{j=0}^{\infty }\sum_{l=0}^{\infty }\int_{0}^{t}s^{\nu -1}%
\frac{(-\lambda s^{\nu })^{j}}{\Gamma (\nu j+\nu )}(t-s)^{\nu -1}\frac{%
(-\lambda (t-s)^{\nu })^{l}}{\Gamma (\nu l+\nu )}ds \\
&=&\lambda ^{2}t^{2\nu -1}\sum_{j=0}^{\infty }\sum_{l=0}^{\infty }\frac{%
(-\lambda t^{\nu })^{j+l}}{\Gamma (\nu j+\nu l+2\nu )} \\
&=&\lambda ^{2}t^{2\nu -1}\sum_{m=0}^{\infty }\frac{(-\lambda t^{\nu })^{m}}{%
\Gamma (\nu m+2\nu )} \\
&=&\lambda ^{2}t^{2\nu -1}E_{\nu ,2\nu }^{2}(-\lambda t^{\nu })=\widehat{f}%
_{1}^{\nu }(t)
\end{eqnarray*}%
\textit{This confirms that one out of two Poisson events are disregarded in
this case (as described in Remark 3.1).}

\

\noindent \textbf{Remark 3.4}\textit{\textbf{\ }We evaluate now the
asymptotic behavior of the interarrival-times density (\ref{lap8}), as
follows:}%
\begin{eqnarray}
\Pr \left\{ \widehat{\mathcal{U}}_{j}\in dt\right\} /dt &=&\lambda
^{2}t^{2\nu -1}\sum_{j=0}^{\infty }\frac{(-\lambda t^{\nu })^{j}(j+1)!}{%
j!\Gamma (\nu j+2\nu )}  \label{rai2} \\
&=&\lambda ^{2}t^{2\nu -1}\sum_{l=0}^{\infty }\frac{l(-\lambda t^{\nu
})^{l-1}}{\Gamma (\nu l+\nu )}=-\frac{\lambda t^{\nu }}{\nu }\frac{d}{dt}%
E_{\nu ,\nu }(-\lambda t^{\nu })  \notag \\
&=&\frac{1-\nu }{\nu }\frac{d}{dt}E_{\nu ,1}(-\lambda t^{\nu })+\frac{t}{\nu
}\frac{d^{2}}{dt^{2}}E_{\nu ,1}(-\lambda t^{\nu }).  \notag
\end{eqnarray}%
\textit{By applying (\ref{rai}) and (\ref{rai3}), we finally get}%
\begin{eqnarray}
\Pr \left\{ \widehat{\mathcal{U}}_{j}\in dt\right\} /dt &=&\frac{\nu -1}{\nu
}\lambda ^{1/\nu }\frac{\sin \left( \nu \pi \right) }{\pi }\int_{0}^{+\infty
}\frac{r^{\nu }e^{-\lambda ^{1/\nu }rt}}{r^{2\nu }+2r^{\nu }\cos (\nu \pi )+1%
}dr+  \notag \\
&&+\frac{\lambda ^{2/\nu }t}{\nu }\frac{\sin \left( \nu \pi \right) }{\pi }%
\int_{0}^{+\infty }\frac{r^{\nu +1}e^{-\lambda ^{1/\nu }rt}}{r^{2\nu
}+2r^{\nu }\cos (\nu \pi )+1}dr  \label{rai4} \\
&\sim &\frac{\sin \left( \nu \pi \right) }{\pi }\frac{2\Gamma (\nu +1)}{%
\lambda t^{\nu +1}}=\frac{2\nu }{\lambda \Gamma (1-\nu )t^{\nu +1}}.  \notag
\end{eqnarray}

\textit{If we compare (\ref{rai4}) with the analogous result (\ref{rai})
obtained for the first model, we can conclude that the interarrival-times
density displays the same asymptotic behavior, with the power law decay (\ref%
{rai4}), so that again the mean waiting time is infinite.}

\textit{For }$t\rightarrow 0$\textit{, instead of (\ref{rairai}), we get in
this case}%
\begin{equation*}
\Pr \left\{ \widehat{\mathcal{U}}_{j}\in dt\right\} /dt\sim \frac{\lambda
^{2}t^{2\nu -1}}{\Gamma (2\nu )}.
\end{equation*}%
\textit{The behavior near the origin of the density of the interarrival
times }$\widehat{\mathcal{U}}_{j}$\textit{\ has a different structure for }$%
\nu <1/2$\textit{\ (tends to infinity) and }$\nu \in \left( 1/2,1\right] $%
\textit{\ (vanishes for }$t\rightarrow 0^{+}$\textit{). For }$\nu =\frac{1}{2%
}$\textit{\ we have instead that }$\Pr \left\{ \widehat{\mathcal{U}}_{j}\in
dt\right\} $\textit{\ is constant for }$t=0$\textit{\ and equal to }$\lambda
^{2}.$

\

We can conclude that $\widehat{\mathcal{N}}_{\nu }$ is a renewal process and
the corresponding renewal function is given by%
\begin{eqnarray}
\widehat{m}_{\nu }(t) &=&\mathbb{E}\widehat{\mathcal{N}}_{\nu }(t)=\lambda
^{2}t^{2\nu }E_{\nu ,2\nu +1}(-2\lambda t^{\nu })  \label{lap9} \\
&=&\frac{\lambda t^{\nu }}{2\Gamma \left( \nu +1\right) }-\frac{\lambda
t^{\nu }}{2}E_{\nu ,\nu +1}(-2\lambda t^{\nu }).  \notag
\end{eqnarray}

Formula (\ref{lap9}) can be obtained either by deriving the probability
generating function (\ref{sec6}) (for $u=1$) or by using the well-known
relationship between the Laplace transforms of the renewal function and the
interarrival-times density:%
\begin{equation*}
\mathcal{L}\left\{ m(t);s\right\} =\frac{1}{s}\frac{\mathcal{L}\left\{
f_{1}(t);s\right\} }{1-\mathcal{L}\left\{ f_{1}(t);s\right\} }.
\end{equation*}%
Indeed, in this case, it is%
\begin{eqnarray*}
\mathcal{L}\left\{ \widehat{m}_{\nu }(t);s\right\} &=&\frac{1}{s}\frac{%
\mathcal{L}\left\{ \widehat{f}_{1}^{\nu }(t);s\right\} }{1-\mathcal{L}%
\left\{ \widehat{f}_{1}^{\nu }(t);s\right\} } \\
&=&\frac{1}{s}\frac{\frac{\lambda ^{2}}{\left( s^{\nu }+\lambda \right) ^{2}}%
}{1-\frac{\lambda ^{2}}{\left( s^{\nu }+\lambda \right) ^{2}}}=\frac{\lambda
^{2}s^{-\nu -1}}{s^{\nu }+2\lambda },
\end{eqnarray*}%
which gives (\ref{lap9}), by applying (\ref{pra}) for $\beta =\nu ,$ $\delta
=1$ and $\gamma =2\nu +1.$

\

\noindent \textbf{Remark 3.5}\textit{\ By comparing the second form of (\ref%
{lap9}) with (\ref{ren}), we can note that the following relationship
between the renewal functions of the two models holds:}%
\begin{equation}
\widehat{m}_{\nu }(t)=\frac{m_{\nu }(t)}{2}-\frac{\lambda t^{\nu }}{2}E_{\nu
,\nu +1}(-2\lambda t^{\nu }).  \label{rn4}
\end{equation}%
\textit{This can be alternatively proved by applying (\ref{sec4}) as follows:%
}%
\begin{eqnarray}
\widehat{m}_{\nu }(t) &=&\sum_{k=0}^{\infty }k\widehat{p}_{k}^{\nu
}(t)=\sum_{k=0}^{\infty }kp_{2k}^{\nu }(t)+\sum_{k=0}^{\infty
}kp_{2k+1}^{\nu }(t)  \label{sec8} \\
&=&\frac{1}{2}\sum_{k=0}^{\infty }(2k)p_{2k}^{\nu }(t)+\frac{1}{2}%
\sum_{k=0}^{\infty }(2k+1)p_{2k+1}^{\nu }(t)-\frac{1}{2}\sum_{k=0}^{\infty
}p_{2k+1}^{\nu }(t)  \notag \\
&=&\frac{1}{2}\sum_{j=0}^{\infty }jp_{j}^{\nu }(t)-\frac{1}{2}%
\sum_{k=0}^{\infty }p_{2k+1}^{\nu }(t)  \notag \\
&=&\frac{m_{\nu }(t)}{2}-\frac{1}{2}\sum_{k=0}^{\infty }p_{2k+1}^{\nu }(t).
\notag
\end{eqnarray}%
\textit{The last term in (\ref{sec8}), which coincides with the sum of the
probabilities of an odd number of events of the first model, can be
evaluated as follows:}%
\begin{eqnarray*}
\sum_{k=0}^{\infty }p_{2k+1}^{\nu }(t) &=&\sum_{k=0}^{\infty }\lambda
^{2k+1}t^{(2k+1)\nu }E_{\nu ,\nu (2k+1)+1}^{2k+2}(-\lambda t^{\nu }) \\
&=&\lambda t^{\nu }E_{\nu ,\nu +1}(-2\lambda t^{\nu }),
\end{eqnarray*}%
\textit{as can be checked by resorting to the Laplace transform:}%
\begin{eqnarray*}
\sum_{k=0}^{\infty }\frac{\lambda ^{2k+1}s^{\nu -1}}{(s^{\nu }+\lambda
)^{2k+2}} &=&\frac{\lambda }{s(s^{\nu }+2\lambda )} \\
&=&\mathcal{L}\left\{ \lambda t^{\nu }E_{\nu ,\nu +1}(-2\lambda t^{\nu
});s\right\} .
\end{eqnarray*}

\textit{Formula (\ref{rn4}) confirms that the mean waiting time is infinite
(as we have noticed in Remark 3.4): indeed it is}%
\begin{equation*}
\lim_{t\rightarrow \infty }\frac{\widehat{m}_{\nu }(t)}{t}%
=\lim_{t\rightarrow \infty }\frac{m_{\nu }(t)}{2t}-\lim_{t\rightarrow \infty
}\frac{\lambda t^{\nu -1}}{2}E_{\nu ,\nu +1}(-2\lambda t^{\nu })=0,
\end{equation*}%
\textit{where the second limit can be evaluated by the following
considerations:}%
\begin{eqnarray*}
\frac{\lambda t^{\nu -1}}{2}E_{\nu ,\nu +1}(-2\lambda t^{\nu }) &=&\frac{1}{%
4t}\left[ 1-E_{\nu ,1}(-2\lambda t^{\nu })\right] \\
&\sim &\frac{1}{4t}\left[ 1-\frac{\sin (\pi \nu )}{\pi }\frac{\Gamma (\nu )}{%
2\lambda t^{\nu }}\right]
\end{eqnarray*}%
\textit{(see (\ref{cau2}) in the Appendix).}

\subsection{The special case $\protect\nu =1$}

We consider now the previous results in the special case $\nu =1$.\ Equation
(\ref{sec1}) reduces in this case to the second-order equation:%
\begin{equation}
\frac{d^{2}\widehat{p}_{k}}{dt^{2}}+2\lambda \frac{d\widehat{p}_{k}}{dt}%
=-\lambda ^{2}(\widehat{p}_{k}-\widehat{p}_{k-1}),\quad k\geq 0  \label{sp1}
\end{equation}%
and the corresponding solution (\ref{sec3}) is given by
\begin{equation}
\widehat{p}_{k}(t)=\frac{(\lambda t)^{2k}}{(2k)!}e^{-\lambda t}+\frac{%
(\lambda t)^{2k+1}}{(2k+1)!}e^{-\lambda t},\quad k\geq 0.  \label{sp2}
\end{equation}%
It is easy to check can show that (\ref{sp2}) solves (\ref{sp1}) with
initial condition%
\begin{eqnarray}
\widehat{p}_{k}(0) &=&\left\{
\begin{array}{c}
1\qquad k=0 \\
0\qquad k\geq 1%
\end{array}%
\right. \\
\widehat{p}_{k}^{\prime }(0) &=&0,\qquad k\geq 0  \notag
\end{eqnarray}%
and $\widehat{p}_{-1}(t)=0.$ Indeed it is%
\begin{eqnarray*}
\frac{d\widehat{p}_{k}}{dt} &=&\frac{\lambda ^{2k}t^{2k-1}}{(2k-1)!}%
e^{-\lambda t}-\frac{\lambda ^{2k+2}t^{2k+1}}{(2k+1)!}e^{-\lambda t} \\
\frac{d^{2}\widehat{p}_{k}}{dt^{2}} &=&\frac{\lambda ^{2k}t^{2k-2}}{(2k-2)!}%
e^{-\lambda t}-\frac{\lambda ^{2k+1}t^{2k-1}}{(2k-1)!}e^{-\lambda t}-\frac{%
\lambda ^{2k+2}t^{2k}}{(2k)!}e^{-\lambda t}+\frac{\lambda ^{2k+3}t^{2k+1}}{%
(2k+1)!}e^{-\lambda t},
\end{eqnarray*}%
so that we get%
\begin{eqnarray*}
&&\frac{d^{2}\widehat{p}_{k}}{dt^{2}}+2\lambda \frac{d\widehat{p}_{k}}{dt} \\
&=&\frac{\lambda ^{2k}t^{2k-2}}{(2k-2)!}e^{-\lambda t}-\frac{\lambda
^{2k+1}t^{2k-1}}{(2k-1)!}e^{-\lambda t}-\frac{\lambda ^{2k+2}t^{2k}}{(2k)!}%
e^{-\lambda t}+\frac{\lambda ^{2k+3}t^{2k+1}}{(2k+1)!}e^{-\lambda t}+ \\
&&+\frac{2\lambda ^{2k+1}t^{2k-1}}{(2k-1)!}e^{-\lambda t}-\frac{2\lambda
^{2k+3}t^{2k+1}}{(2k+1)!}e^{-\lambda t} \\
&=&\frac{\lambda ^{2k}t^{2k-2}}{(2k-2)!}e^{-\lambda t}+\frac{\lambda
^{2k+1}t^{2k-1}}{(2k-1)!}e^{-\lambda t}-\frac{\lambda ^{2k+2}t^{2k}}{(2k)!}%
e^{-\lambda t}-\frac{\lambda ^{2k+3}t^{2k+1}}{(2k+1)!}e^{-\lambda t} \\
&=&-\lambda ^{2}(\widehat{p}_{k}-\widehat{p}_{k-1}).
\end{eqnarray*}

Formula (\ref{sp2}) agrees with the relationship (\ref{sec4}) given in
Remark 3.1, which in this case we can write as follows%
\begin{equation}
\widehat{p}_{k}(t)=p_{2k}(t)+p_{2k+1}(t).  \label{rel8}
\end{equation}%
Therefore it can be interpreted, for $k=0,1,2,...,$\ as the probability
distribution $\Pr \left\{ \widehat{N}(t)=k\right\} $\ for a
\textquotedblleft second-order process\textquotedblright\ $\widehat{N}%
(t),t>0\ $linked to the standard Poisson process $N$ by the following
relationship:%
\begin{equation}
\Pr \left\{ \widehat{N}(t)=k\right\} =\Pr \left\{ N(t)=2k\right\} +\Pr
\left\{ N(t)=2k+1\right\} .  \label{rel9}
\end{equation}%
From (\ref{lap6}) and (\ref{lap8}) we can easily see that for this process
the densities of the interarrival times and of the $k$-th event waiting time
are given respectively by%
\begin{equation}
\widehat{f}_{1}(t)=\lambda ^{2}te^{-\lambda t}  \label{sp3}
\end{equation}%
and
\begin{equation}
\widehat{f}_{k}(t)=\frac{\lambda ^{2k}t^{2k-1}}{(2k-1)!}e^{-\lambda t}.
\label{sp4}
\end{equation}%
Therefore, in this case, the random variable $\widehat{T}_{j},$ representing
the instant of the $j$-th event, is distributed as $Gamma(\lambda ,2j).$

We derive equation (\ref{sp2}) in an alternative way, which is similar to
the construction of the standard Poisson process, by considering the
relationships (\ref{sec9}) or, equivalently, the property of the
interarrival times described in Remark 3.3. We can write that

\begin{equation}
\widehat{T}_{j}\overset{law}{=}T_{2j},  \label{rel10}
\end{equation}%
where%
\begin{equation*}
T_{j}=\inf \left\{ t>0:N(t)=j\right\}
\end{equation*}%
which represents the time of the $j$-th event of $N(t).$

.

\noindent \textbf{Theorem 3.3 \ }The probability distribution of the process
$\widehat{N}(t),t>0$ described by (\ref{rel9}) and (\ref{rel10}) solves
equation (\ref{sp1}).

\noindent \textbf{Proof \ }Let us consider now the following intervals

\begin{eqnarray*}
A &=&\left[ T_{2k-2},T_{2k-1}\right) \\
B &=&\left[ T_{2k-1},T_{2k}\right) \\
C &=&\left[ T_{2k},T_{2k+1}\right) \\
D &=&\left[ T_{2k+1},T_{2k+2}\right)
\end{eqnarray*}%
so that $\left[ T_{2k-2},T_{2k+2}\right) =A\cup B\cup C\cup D$ (see Fig.1).

\begin{figure}[h]

\centering {\includegraphics[scale=0.5]{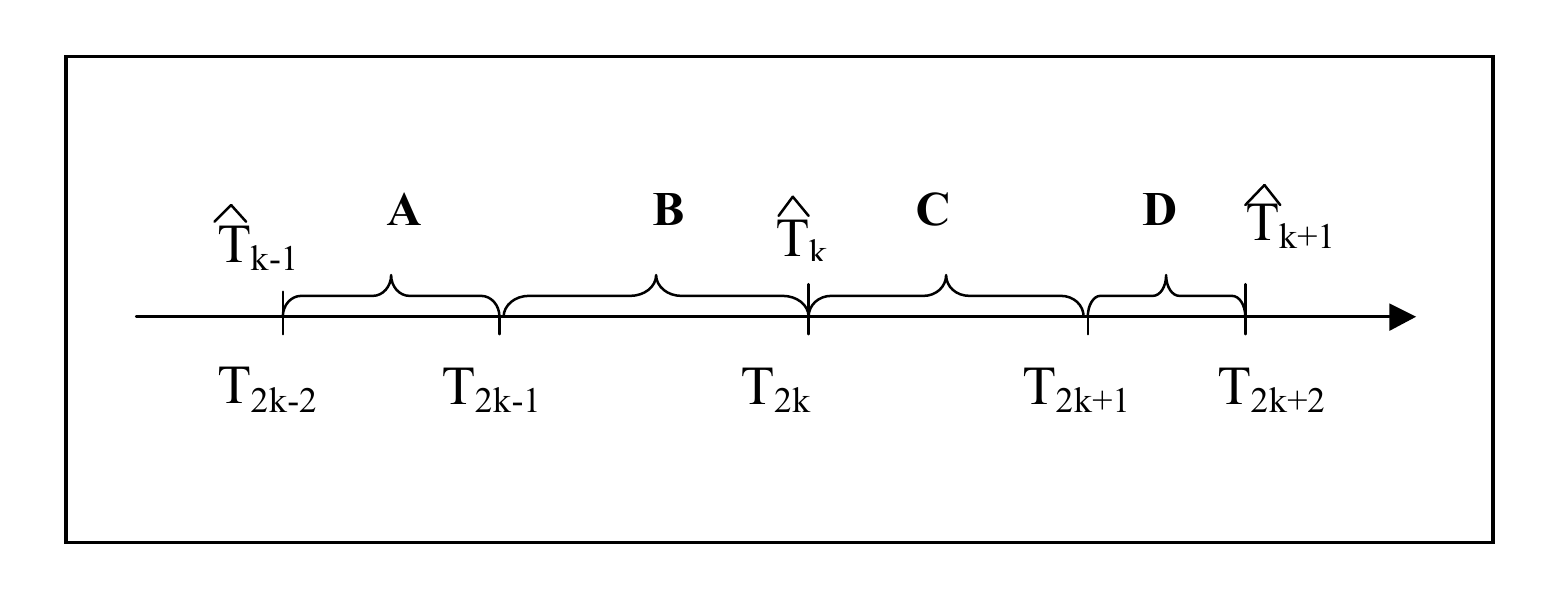}}
\caption{The interval $\left[ T_{2k-2},T_{2k+2}\right)$}

\end{figure}

We evaluate the following probability, by stopping the approximation at the
second-order terms:

\begin{eqnarray}
&&\widehat{p}_{k}(t+2\Delta t)  \label{app} \\
&\simeq &2\Pr \left( \widehat{N}(t)=k-1,t\in A\right) \frac{\lambda
^{2}\left( \Delta t\right) ^{2}}{2}\left( 1-\frac{\lambda ^{2}\left( \Delta
t\right) ^{2}}{2}\right) +  \notag \\
&&+\Pr \left( \widehat{N}(t)=k-1,t\in A\right) \left( \lambda \Delta
t-\lambda ^{2}\left( \Delta t\right) ^{2}\right) ^{2}+  \notag \\
&&+2\Pr \left( \widehat{N}(t)=k-1,t\in B\right) \left( \lambda \Delta
t-\lambda ^{2}\left( \Delta t\right) ^{2}\right) \left( 1-\lambda \Delta t+%
\frac{\lambda ^{2}\left( \Delta t\right) ^{2}}{2}\right) +  \notag \\
&&+2\Pr \left( \widehat{N}(t)=k-1,t\in B\right) \frac{\lambda ^{2}\left(
\Delta t\right) ^{2}}{2}\left( 1-\lambda \Delta t+\frac{\lambda ^{2}\left(
\Delta t\right) ^{2}}{2}\right) +  \notag
\end{eqnarray}%
\begin{eqnarray}
&&+\Pr \left( \widehat{N}(t)=k-1,t\in B\right) \left( \lambda \Delta
t-\lambda ^{2}\left( \Delta t\right) ^{2}\right) ^{2}+  \notag \\
&&+\Pr \left( \widehat{N}(t)=k,t\in C\right) \left( 1-\lambda \Delta t+\frac{%
\lambda ^{2}\left( \Delta t\right) ^{2}}{2}\right) ^{2}+  \notag \\
&&+2\Pr \left( \widehat{N}(t)=k,t\in C\right) \left( 1-\lambda \Delta t+%
\frac{\lambda ^{2}\left( \Delta t\right) ^{2}}{2}\right) \left( \lambda
\Delta t-\lambda ^{2}\left( \Delta t\right) ^{2}\right) +  \notag \\
&&+\Pr \left( \widehat{N}(t)=k,t\in D\right) \left( 1-\lambda \Delta t+\frac{%
\lambda ^{2}\left( \Delta t\right) ^{2}}{2}\right) ^{2}+o\left( (\Delta
t)^{2}\right) ,  \notag
\end{eqnarray}%
where we used the following well-known approximations valid for the standard
Poisson process:%
\begin{eqnarray*}
\Pr \left( 0\text{ Poisson events in }\Delta t\right) &=&e^{-\lambda \Delta
t}\simeq 1-\lambda \Delta t+\frac{\lambda ^{2}\left( \Delta t\right) ^{2}}{2}%
+o\left( (\Delta t)^{2}\right) \\
\Pr \left( 1\text{ Poisson event in }\Delta t\right) &=&\lambda \Delta
t\,e^{-\lambda \Delta t}\simeq \lambda \Delta t-\lambda ^{2}\left( \Delta
t\right) ^{2}+o\left( (\Delta t)^{2}\right) \\
\Pr \left( 2\text{ Poisson events in }\Delta t\right) &=&\frac{\lambda
^{2}\left( \Delta t\right) ^{2}}{2}e^{-\lambda \Delta t}\simeq \frac{\lambda
^{2}\left( \Delta t\right) ^{2}}{2}+o\left( (\Delta t)^{2}\right) .
\end{eqnarray*}%
By ignoring the terms of order greater than $(\Delta t)^{2}$ the probability
(\ref{app}) can be rewritten as follows%
\begin{eqnarray}
&&\widehat{p}_{k}(t+2\Delta t)  \label{app2} \\
&\simeq &2\Pr \left( \widehat{N}(t)=k-1,t\in A\right) \frac{\lambda
^{2}\left( \Delta t\right) ^{2}}{2}+  \notag \\
&&+\Pr \left( \widehat{N}(t)=k-1,t\in A\right) \lambda ^{2}\left( \Delta
t\right) ^{2}+  \notag \\
&&+2\Pr \left( \widehat{N}(t)=k-1,t\in B\right) \left( \lambda \Delta
t-2\lambda ^{2}\left( \Delta t\right) ^{2}\right) +  \notag \\
&&+2\Pr \left( \widehat{N}(t)=k-1,t\in B\right) \frac{\lambda ^{2}\left(
\Delta t\right) ^{2}}{2}+  \notag
\end{eqnarray}%
\begin{eqnarray}
&&+\Pr \left( \widehat{N}(t)=k-1,t\in B\right) \lambda ^{2}\left( \Delta
t\right) ^{2}+  \notag \\
&&+\Pr \left( \widehat{N}(t)=k,t\in C\right) \left( 1-2\lambda \Delta
t+2\lambda ^{2}\left( \Delta t\right) ^{2}\right) +  \notag \\
&&+2\Pr \left( \widehat{N}(t)=k,t\in C\right) \left( \lambda \Delta
t-2\lambda ^{2}\left( \Delta t\right) ^{2}\right) +  \notag \\
&&+\Pr \left( \widehat{N}(t)=k,t\in D\right) \left( 1+2\lambda ^{2}\left(
\Delta t\right) ^{2}-2\lambda \Delta t\right) ^{2}+o\left( (\Delta
t)^{2}\right) .  \notag
\end{eqnarray}%
Since we can write%
\begin{equation*}
\Pr \left( \widehat{N}(t)=k-1,t\in A\right) +\Pr (\widehat{N}(t)=k-1,t\in B)=%
\widehat{p}_{k-1}(t)
\end{equation*}%
and, analogously%
\begin{equation*}
\Pr \left( \widehat{N}(t)=k,t\in C\right) +\Pr (\widehat{N}(t)=k,t\in D)=%
\widehat{p}_{k}(t)
\end{equation*}%
formula (\ref{app2}) becomes%
\begin{eqnarray}
\widehat{p}_{k}(t+2\Delta t) &\simeq &2\widehat{p}_{k-1}(t)\frac{\lambda
^{2}\left( \Delta t\right) ^{2}}{2}+\widehat{p}_{k-1}(t)\lambda ^{2}\left(
\Delta t\right) ^{2}+  \label{app4} \\
&+&\widehat{p}_{k}(t)\left( 1-\lambda \Delta t\right) ^{2}+\widehat{p}%
_{k}(t)\lambda ^{2}\left( \Delta t\right) ^{2}+  \notag \\
&&+2\Pr \left( \widehat{N}(t)=k-1,t\in B\right) \left( \lambda \Delta
t-2\lambda ^{2}\left( \Delta t\right) ^{2}\right) +  \notag \\
&&+2\Pr \left( \widehat{N}(t)=k,t\in C\right) \left( \lambda \Delta
t-2\lambda ^{2}\left( \Delta t\right) ^{2}\right) +o\left( (\Delta
t)^{2}\right) .  \notag
\end{eqnarray}%
We consider now the following probability, on a single interval of length $%
\Delta t$
\begin{eqnarray}
&&\widehat{p}_{k}(t+\Delta t)\simeq \Pr \left( \widehat{N}(t)=k-1,t\in
A\right) \frac{\lambda ^{2}\left( \Delta t\right) ^{2}}{2}+  \label{app3} \\
&+&\Pr \left( \widehat{N}(t)=k-1,t\in B\right) \left( \lambda \Delta t-\frac{%
\lambda ^{2}\left( \Delta t\right) ^{2}}{2}\right) +  \notag \\
&+&\Pr \left( \widehat{N}(t)=k,t\in C\right) \left( 1-\frac{\lambda
^{2}\left( \Delta t\right) ^{2}}{2}\right) +  \notag \\
&+&\Pr \left( \widehat{N}(t)=k,t\in D\right) \left( 1-\lambda \Delta t+\frac{%
\lambda ^{2}\Delta t^{2}}{2}\right) +o\left( (\Delta t)^{2}\right)  \notag \\
&=&\left[ \widehat{p}_{k-1}(t)-\Pr \left( \widehat{N}(t)=k-1,t\in B\right) %
\right] \frac{\lambda ^{2}\left( \Delta t\right) ^{2}}{2}+  \notag
\end{eqnarray}%
\begin{eqnarray}
&+&\Pr \left( \widehat{N}(t)=k-1,t\in B\right) \left( \lambda \Delta t-\frac{%
\lambda ^{2}\left( \Delta t\right) ^{2}}{2}\right) +  \notag \\
&+&\Pr \left( \widehat{N}(t)=k,t\in C\right) \left( 1-\frac{\lambda
^{2}\left( \Delta t\right) ^{2}}{2}\right) +  \notag \\
&+&\left[ \widehat{p}_{k}(t)-\Pr \left( \widehat{N}(t)=k,t\in C\right) %
\right] \left( 1-\lambda \Delta t+\frac{\lambda ^{2}\Delta t^{2}}{2}\right)
+o\left( (\Delta t)^{2}\right)  \notag \\
&=&\widehat{p}_{k-1}(t)\frac{\lambda ^{2}\left( \Delta t\right) ^{2}}{2}+\Pr
\left( \widehat{N}(t)=k-1,t\in B\right) \left( \lambda \Delta t-\lambda
^{2}\left( \Delta t\right) ^{2}\right)  \notag \\
&+&\widehat{p}_{k}(t)\left( 1-\lambda \Delta t+\frac{\lambda ^{2}\Delta t^{2}%
}{2}\right) +\Pr \left( \widehat{N}(t)=k,t\in C\right) \left( \lambda \Delta
t-\lambda ^{2}\left( \Delta t\right) ^{2}\right) +o\left( \left( \Delta
t\right) ^{2}\right) .  \notag
\end{eqnarray}%
We multiply (\ref{app3}) by $2(1-\lambda \Delta t)$ and ignore the terms of
order greater than $\left( \Delta t\right) ^{2}$, so that we get%
\begin{eqnarray*}
&&\widehat{p}_{k-1}(t)\lambda ^{2}\left( \Delta t\right) ^{2}\simeq
2(1-\lambda \Delta t)\widehat{p}_{k}(t+\Delta t)-2\Pr \left( \widehat{N}%
(t)=k-1,t\in B\right) \left( \lambda \Delta t-2\lambda ^{2}\left( \Delta
t\right) ^{2}\right) + \\
&-&2\widehat{p}_{k}(t)\left[ \left( 1-\lambda \Delta t\right) ^{2}+\frac{%
\lambda ^{2}\Delta t^{2}}{2}\right] -2\Pr \left( \widehat{N}(t)=k,t\in
C\right) \left( \lambda \Delta t-2\lambda ^{2}\left( \Delta t\right)
^{2}\right) +o\left( \left( \Delta t\right) ^{2}\right) ,
\end{eqnarray*}%
which can be substituted in the first line of (\ref{app4}). Finally we obtain%
\begin{eqnarray}
&&\widehat{p}_{k}(t+2\Delta t)\simeq \widehat{p}_{k-1}(t)\lambda ^{2}\left(
\Delta t\right) ^{2}+2(1-\lambda \Delta t)\widehat{p}_{k}(t+\Delta t)+
\label{app5} \\
&-&2\Pr \left( \widehat{N}(t)=k-1,t\in B\right) \left( \lambda \Delta
t-2\lambda ^{2}\left( \Delta t\right) ^{2}\right) -2\widehat{p}_{k}(t)\left[
\left( 1-\lambda \Delta t\right) ^{2}+\frac{\lambda ^{2}\Delta t^{2}}{2}%
\right] +  \notag \\
&-&2\Pr \left( \widehat{N}(t)=k,t\in C\right) \left( \lambda \Delta
t-2\lambda ^{2}\left( \Delta t\right) ^{2}\right) +  \notag
\end{eqnarray}%
\begin{eqnarray}
&+&2\Pr \left( \widehat{N}(t)=k-1,t\in B\right) \left( \lambda \Delta
t-2\lambda ^{2}\left( \Delta t\right) ^{2}\right) +  \notag \\
&+&2\Pr \left( \widehat{N}(t)=k,t\in C\right) \left( \lambda \Delta
t-2\lambda ^{2}\left( \Delta t\right) ^{2}\right) +  \notag \\
&+&\widehat{p}_{k}(t)\left( 1-\lambda \Delta t\right) ^{2}+\widehat{p}%
_{k}(t)\lambda ^{2}\left( \Delta t\right) ^{2}+o\left( (\Delta t)^{2}\right)
\notag \\
&=&\widehat{p}_{k-1}(t)\lambda ^{2}\left( \Delta t\right) ^{2}+2(1-\lambda
\Delta t)\widehat{p}_{k}(t+\Delta t)-\widehat{p}_{k}(t)\left( 1-\lambda
\Delta t\right) ^{2}.  \notag
\end{eqnarray}%
We divide (\ref{app5}) by $(\Delta t)^{2}$, so that we get%
\begin{equation*}
\frac{\widehat{p}_{k}(t+2\Delta t)-2\widehat{p}_{k}(t+\Delta t)+\widehat{p}%
_{k}(t)}{(\Delta t)^{2}}+2\lambda \Delta t\frac{\widehat{p}_{k}(t+\Delta t)-%
\widehat{p}_{k}(t)}{(\Delta t)^{2}}=-\lambda ^{2}\left( \Delta t\right) ^{2}%
\frac{\widehat{p}_{k}(t)-\widehat{p}_{k-1}(t)}{(\Delta t)^{2}}.
\end{equation*}%
By letting $\Delta t\rightarrow 0$ we easily obtain equation (\ref{sp1})$.$%
\hfill $\square $

\

\noindent \textbf{Remark 3.6}\textit{\ This special case is particularly
interesting because it describes a generalization of the Poisson process
which have been used in many articles. Random motions at finite velocities
spaced by this particular renewal process have been considered by different
authors ([5]-[6]-[10]).}

\textit{In particular in [3] it has been studied a model with uniformly
distributed deviations which take place at even-order Poisson events and
therefore its interarrival times are distributed as }$Gamma(\lambda ,2).$%
\textit{\ }

\section{Conclusions}

The results of the previous sections can be generalized to the $n$-th order
case, if we consider the following equation%
\begin{equation}
\frac{d^{n\nu }p_{k}}{dt^{n\nu }}+\binom{n}{1}\lambda \frac{d^{(n-1)\nu
}p_{k}}{dt^{(n-1)\nu }}+...+\binom{n}{n-1}\lambda ^{n-1}\frac{d^{\nu }p_{k}}{%
dt^{\nu }}=-\lambda ^{n}(p_{k}-p_{k-1}),\quad k\geq 0,  \label{n}
\end{equation}%
where $\nu \in \left( 0,1\right) ,$ subject to the initial conditions%
\begin{eqnarray}
p_{k}(0) &=&\left\{
\begin{array}{c}
1\qquad k=0 \\
0\qquad k\geq 1%
\end{array}%
\right. ,\text{ \quad for }0<\nu <1  \label{n1} \\
\left. \frac{d^{j}}{dt^{j}}p_{k}(t)\right\vert _{t=0} &=&0\qquad
j=1,...,n-1,\quad k\geq 0,\text{ \quad for }\frac{1}{n}<\nu <1  \notag
\end{eqnarray}%
and $p_{-1}(t)=0$. Following the same steps as for the first two models we
get the Laplace transform of (\ref{n}):%
\begin{equation*}
\left[ s^{n\nu }+\binom{n}{1}\lambda s^{(n-1)\nu }+...+\binom{n}{n-1}\lambda
^{n-1}s^{2\nu }+\lambda ^{n}\right] \mathcal{L}\left\{ \widetilde{p}%
_{k}^{\nu }(t);s\right\} =\lambda ^{n}\mathcal{L}\left\{ \widetilde{p}%
_{k-1}^{\nu }(t);s\right\} ,
\end{equation*}%
which can be solved, in view of (\ref{jen}) and the initial conditions (\ref%
{n1}), recursively, yielding%
\begin{equation}
\mathcal{L}\left\{ \widetilde{p}_{k}^{\nu }(t);s\right\} =\frac{%
\sum_{j=1}^{n}\binom{n}{j}s^{\nu j-1}\lambda ^{(k+1)n-j}}{(s^{\nu }+\lambda
)^{(k+1)n}}.  \label{n2}
\end{equation}%
For $n=2,$ we obtain from (\ref{n2}) formula (\ref{lap}). The Laplace
transform can be inverted by applying again (\ref{pra}) and the solution is
given, also in this case, as a sum of GML functions as follows:%
\begin{equation}
\widetilde{p}_{k}^{\nu }(t)=\sum_{j=1}^{n}\binom{n}{j}\left( \lambda t^{\nu
}\right) ^{n(k+1)-j}E_{\nu ,\nu n(k+1)-\nu j+1}^{kn+n}(-\lambda t^{\nu }).
\label{p}
\end{equation}%
For $n=1$, we get the distribution of the first model (\ref{gml3}), while,
for $n=2$, we get the solution of the second-type equation in the form (\ref%
{pre}). In this case the use of (\ref{gen}) requires much harder
calculations. Nevertheless a relationship similar to (\ref{sec4}) can be
obtained, even for $n>2$, by studying the density $\widetilde{f}_{k}^{\nu
}(t),t>0$ of the waiting time of the $k$-th event $T_{k}$. As already seen
in section 3, the following identity must be satisfied by $\widetilde{f}%
_{k}^{\nu }(t):$
\begin{equation}
\mathcal{L}\left\{ \widetilde{p}_{k}^{\nu }(t);s\right\} =\frac{1}{s}%
\mathcal{L}\left\{ \widetilde{f}_{k}^{\nu }(t);s\right\} -\frac{1}{s}%
\mathcal{L}\left\{ \widetilde{f}_{k+1}^{\nu }(t);s\right\} ,  \label{rel7}
\end{equation}%
so that, by substituting (\ref{n2}) in the l.h.s. of (\ref{rel7}) we get%
\begin{equation}
\mathcal{L}\left\{ \widetilde{f}_{k}^{\nu }(t);s\right\} =\frac{\lambda ^{nk}%
}{(s^{\nu }+\lambda )^{nk}},
\end{equation}%
which can be inverted as usual, thus obtaining%
\begin{equation}
\widetilde{f}_{k}^{\nu }(t)=\lambda ^{nk}t^{n\nu k-1}E_{\nu ,n\nu
k}^{nk}(-\lambda t^{\nu }).  \label{p2}
\end{equation}%
Again the process is a renewal one, since (\ref{p2}) coincides with the sum
of $k$ independent and identically distributed random variables $\widetilde{%
\mathcal{U}}_{j}$'s (representing the interarrival times) with density given
by%
\begin{equation}
\Pr \left\{ \widetilde{\mathcal{U}}_{j}\in dt\right\} /dt=\mathcal{L}%
^{-1}\left\{ \frac{\lambda ^{n}}{\left( s^{\nu }+\lambda \right) ^{n}}%
;t\right\} =\lambda ^{n}t^{n\nu -1}E_{\nu ,n\nu }^{n}(-\lambda t^{\nu
})=f_{1}^{n\nu }(t).  \label{ul}
\end{equation}%
Formula (\ref{ul}) shows that each interarrival time of the $n$-th order
case is distributed as the sum of $n$ independent interarrival times of the
first model. This suggests that the following relationship between the
corresponding probability distributions holds:%
\begin{equation}
\widetilde{p}_{k}^{\nu }(t)=p_{nk}^{\nu }(t)+p_{nk+1}^{\nu
}(t)+...+p_{nk+n-1}^{\nu }(t),\quad n>2.
\end{equation}

\section{Appendix}

We derive the integral form of the Mittag-Leffler function (used in (\ref%
{rai3})), and some generalizations. We start by showing that, for $0<\nu <1$,%
\begin{equation}
E_{\nu ,1}(-t^{\nu })=\frac{\sin \left( \nu \pi \right) }{\pi }%
\int_{0}^{+\infty }\frac{r^{\nu -1}e^{-rt}}{r^{2\nu }+2r^{\nu }\cos (\nu \pi
)+1}dr.  \label{rai5}
\end{equation}%
By applying the reflection formula of the Gamma function%
\begin{equation*}
\Gamma (\nu m+1)\Gamma (-\nu m)=\frac{\pi }{\sin (-\pi \nu m)}
\end{equation*}%
we can write the Mittag-Leffler function as follows%
\begin{eqnarray*}
E_{\nu ,1}(-t^{\nu }) &=&-\frac{1}{\pi }\sum_{m=0}^{\infty }\frac{(-t^{\nu
})^{m}\sin (\pi \nu m)}{\Gamma (\nu m+1)}\Gamma (\nu m+1)\Gamma (-\nu m) \\
&=&-\frac{1}{\pi }\sum_{m=0}^{\infty }\frac{(-1)^{m}\sin (\pi \nu m)}{\Gamma
(\nu m+1)}\int_{0}^{+\infty }e^{-rt}r^{-\nu m-1}dr\int_{0}^{+\infty
}e^{-y}y^{\nu m}dy \\
&=&-\frac{1}{\pi }\sum_{m=0}^{\infty }\frac{(-1)^{m}\sin (\pi \nu m)}{\Gamma
(\nu m+1)}\int_{0}^{+\infty }e^{-rt}\left( \int_{0}^{+\infty }e^{-ry}y^{\nu
m}dy\right) dr \\
&=&-\frac{1}{\pi }\int_{0}^{+\infty }e^{-rt}\left( \int_{0}^{+\infty
}e^{-ry}\sum_{m=0}^{\infty }\frac{(-1)^{m}y^{\nu m}}{\Gamma (\nu m+1)}\frac{%
e^{i\pi \nu m}-e^{-i\pi \nu m}}{2i}dy\right) dr \\
&=&-\frac{1}{2\pi i}\int_{0}^{+\infty }e^{-rt}\left( \int_{0}^{+\infty
}e^{-ry}E_{\nu ,1}(-y^{\nu }e^{i\pi \nu })dy-\int_{0}^{+\infty
}e^{-ry}E_{\nu ,1}(-y^{\nu }e^{-i\pi \nu })dy\right) dr \\
&=&-\frac{1}{2\pi i}\int_{0}^{+\infty }r^{\nu -1}e^{-rt}\left[ \frac{1}{%
r^{\nu }+e^{i\pi \nu }}-\frac{1}{r^{\nu }+e^{-i\pi \nu }}\right] dr \\
&=&\frac{\sin \left( \nu \pi \right) }{\pi }\int_{0}^{+\infty }\frac{r^{\nu
-1}e^{-rt}}{(r^{\nu }+e^{i\pi \nu })(r^{\nu }+e^{-i\pi \nu })}dr.
\end{eqnarray*}

From (\ref{rai5}) we can derive the following approximation, for large $t$:%
\begin{equation}
E_{\nu ,1}(-t^{\nu })\simeq \frac{\sin \left( \nu \pi \right) }{\pi }\frac{%
\Gamma (\nu )}{t^{\nu }}.  \label{cau2}
\end{equation}

We present now an expression of the Mittag-Leffler function which permits us
to interpret it as the mean of a Cauchy r.v. We can rewrite (\ref{rai5}) as

\begin{eqnarray}
E_{\nu ,1}(-t^{\nu }) &=&\frac{1}{\pi }\int_{0}^{+\infty }\frac{\sin \left(
\nu \pi \right) }{[r^{\nu }+\cos (\nu \pi )]^{2}+\sin ^{2}(\nu \pi )}r^{\nu
-1}e^{-rt}dr  \label{cau} \\
&=&\frac{1}{\pi \nu }\int_{0}^{+\infty }\frac{\sin \left( \nu \pi \right) }{%
[r+\cos (\nu \pi )]^{2}+\sin ^{2}(\nu \pi )}e^{-r^{1/\nu }t}dr  \notag \\
&=&\mathbb{E}_{X}\left\{ \frac{1}{\nu }e^{-tX^{1/\nu }}1_{\left[ 0,\infty
\right) }\right\} ,  \notag
\end{eqnarray}%
where $X$ is distributed as a $Cauchy$ with parameters $-\cos (\pi \nu )$
and $\sin \left( \nu \pi \right) .$ Formula (\ref{cau}) permits us to study
the particular case where $\nu =1$, since we can write (\ref{cau}), by means
of the characteristic function of a Cauchy r.v., as follows:%
\begin{gather}
E_{\nu ,1}(-t^{\nu })=\frac{1}{2\pi \nu }\int_{0}^{+\infty }e^{-r^{1/\nu
}t}\left( \int_{-\infty }^{+\infty }e^{-ir\beta -|\beta |\sin (\pi \nu
)-i\beta \cos (\pi \nu )}d\beta \right) dr  \label{ni} \\
\overset{\nu \rightarrow 1^{-}}{\rightarrow }\int_{0}^{+\infty
}e^{-rt}\delta (r-1)dr=e^{-t}=E_{1,1}(-t).  \notag
\end{gather}%
Analogously we can study the case where $\nu =1/2$: from the first line of (%
\ref{ni}) and by applying formula 3.896.4, p.480 of [7], we get%
\begin{eqnarray*}
&&E_{\frac{1}{2},1}(-\sqrt{t})=\frac{1}{\pi }\int_{0}^{+\infty
}e^{-r^{2}t}\left( \int_{-\infty }^{+\infty }e^{-ir\beta -|\beta |}d\beta
\right) dr \\
&=&\frac{1}{\pi }\int_{0}^{+\infty }e^{-r^{2}t}\left( \int_{0}^{+\infty
}e^{-ir\beta -\beta }d\beta +\int_{-\infty }^{0}e^{-ir\beta +\beta }d\beta
\right) dr \\
&=&\frac{2}{\pi }\int_{0}^{+\infty }e^{-r^{2}t}\cos (\beta r)\left(
\int_{0}^{+\infty }e^{-\beta }d\beta \right) dr \\
&=&\frac{2}{\pi }\int_{0}^{+\infty }e^{-\beta }\left( \int_{0}^{+\infty
}e^{-r^{2}t}\cos (\beta r)dr\right) d\beta \\
&=&\frac{1}{\pi }\sqrt{\frac{\pi }{t}}\int_{0}^{+\infty }e^{-\beta }e^{-%
\frac{\beta ^{2}}{4t}}d\beta \\
&=&[y^{2}=\frac{\beta ^{2}}{4t}] \\
&=&\frac{2}{\sqrt{\pi }}\int_{0}^{+\infty }e^{-y^{2}-2\sqrt{t}}dy,
\end{eqnarray*}%
which coincides with the form given in [17], for $x=-\sqrt{t}$.

We generalize formula (\ref{rai5}) to the case of a two-parameters
Mittag-Leffler function: for $0<\nu <1$ and $0<\beta <\nu +1$, we have that%
\begin{equation}
E_{\nu ,\beta }(-t^{\nu })=\frac{t^{1-\beta }}{\pi }\int_{0}^{+\infty
}r^{\nu -\beta }e^{-rt}\frac{\sin \left( \nu \pi \right) }{\left[ r^{\nu
}+\cos (\nu \pi )\right] ^{2}+\sin ^{2}(\nu \pi )}\left[ \frac{\sin \left(
\beta \pi \right) }{\sin \left( \nu \pi \right) }\left[ r^{\nu }-\cos (\nu
\pi )\right] +\cos (\beta \pi )\right] dr.  \label{rai8}
\end{equation}

We derive (\ref{rai8}) by starting from the series representation of the
Mittag-Leffler function:

\begin{eqnarray*}
&&E_{\nu ,\beta }(-t^{\nu }) \\
&=&\sum_{m=0}^{\infty }\frac{(-1)^{m}t^{\nu m}}{\Gamma (\nu m+\beta )}\frac{%
\sin (\left( \nu m+\beta )\pi \right) }{\pi }\frac{\pi }{\sin (\left( \nu
m+\beta )\pi \right) } \\
&=&\frac{1}{\pi }\sum_{m=0}^{\infty }\frac{(-1)^{m}t^{\nu m}}{\Gamma (\nu
m+\beta )}\frac{\sin (\left( \nu m+\beta )\pi \right) }{\pi }\Gamma (1-\nu
m-\beta )\Gamma (\nu m+\beta ) \\
&=&\frac{t^{1-\beta }}{\pi }\sum_{m=0}^{\infty }\frac{(-1)^{m}}{\Gamma (\nu
m+\beta )}\sin (\left( \nu m+\beta )\pi \right) \int_{0}^{+\infty
}e^{-rt}r^{-\nu m-\beta }dr\int_{0}^{+\infty }e^{-y}y^{\nu m+\beta -1}dy \\
&=&\frac{t^{1-\beta }}{\pi }\sum_{m=0}^{\infty }\frac{(-1)^{m}}{\Gamma (\nu
m+\beta )}\sin (\left( \nu m+\beta )\pi \right) \int_{0}^{+\infty
}e^{-rt}\left( \int_{0}^{+\infty }e^{-ry}y^{\nu m+\beta -1}dy\right) dr \\
&=&\frac{t^{1-\beta }}{\pi }\int_{0}^{+\infty }e^{-rt}\left(
\int_{0}^{+\infty }e^{-ry}y^{\beta -1}\sum_{m=0}^{\infty }\frac{%
(-1)^{m}y^{\nu m}}{\Gamma (\nu m+\beta )}\frac{e^{i\pi \nu m+i\pi \beta
}-e^{-i\pi \nu m-i\pi \beta }}{2i}dy\right) dr
\end{eqnarray*}%
\begin{eqnarray*}
&=&\frac{t^{1-\beta }}{2\pi i}\int_{0}^{\infty }dy\int_{0}^{\infty
}e^{-r(y+t)}y^{\beta -1}\left[ e^{i\pi \beta }E_{\nu ,\beta }(-y^{\nu
}e^{i\pi \nu })-e^{-i\pi \beta }E_{\nu ,\beta }(-y^{\nu }e^{-i\pi \nu })%
\right] dr \\
&=&\frac{t^{1-\beta }}{2\pi i}\int_{0}^{\infty }e^{-rt}\left[ e^{i\pi \beta }%
\frac{r^{\nu -\beta }}{r^{\nu }+e^{i\pi \nu }}-e^{-i\pi \beta }\frac{r^{\nu
-\beta }}{r^{\nu }+e^{-i\pi \nu }}\right] dr \\
&=&\frac{t^{1-\beta }}{\pi }\int_{0}^{\infty }e^{-rt}r^{\nu -\beta }\frac{%
r^{\nu }\sin (\pi \beta )+\sin (\pi (\beta -\nu ))}{r^{2\nu }+2r^{\nu }\cos
(\pi \nu )+1}dr \\
&=&\frac{t^{1-\beta }}{\pi }\int_{0}^{\infty }e^{-rt}r^{\nu -\beta }\frac{%
\sin (\pi \nu )}{(r^{\nu }+\cos (\pi \nu ))^{2}+\sin ^{2}(\pi \nu )}\left[
r^{\nu }\frac{\sin (\pi \beta )}{\sin (\pi \nu )}+\frac{\sin (\pi \beta
)\cos (\pi \nu )-\sin (\pi \nu )\cos (\pi \beta )}{\sin (\pi \nu )}\right]
dr.
\end{eqnarray*}%
From the previous expression formula (\ref{rai8}) easily follows and, for $%
\beta =1$, it reduces to (\ref{rai5}). For $\beta =\nu $ we obtain from (\ref%
{rai8})%
\begin{equation}
E_{\nu ,\nu }(-t^{\nu })=\frac{t^{1-\nu }}{\pi }\int_{0}^{+\infty }r^{\nu
}e^{-rt}\frac{\sin \left( \nu \pi \right) }{\left[ r^{\nu }+\cos (\nu \pi )%
\right] ^{2}+\sin ^{2}(\nu \pi )}dr.  \label{rai9}
\end{equation}%
As a check of (\ref{rai9}) we can study the limit for $\nu \rightarrow 1:$%
\begin{eqnarray*}
E_{\nu ,\nu }(-t^{\nu }) &=&\frac{t^{1-\nu }}{\nu \pi }\int_{0}^{+\infty }r^{%
\frac{1}{\nu }}e^{-r^{\frac{1}{\nu }}t}\frac{\sin \left( \nu \pi \right) }{%
\left[ r+\cos (\nu \pi )\right] ^{2}+\sin ^{2}(\nu \pi )}dr \\
&=&\frac{t^{1-\nu }}{2\pi \nu }\int_{0}^{+\infty }r^{\frac{1}{\nu }}e^{-r^{%
\frac{1}{\nu }}t}\left( \int_{-\infty }^{+\infty }e^{-ir\beta -|\beta |\sin
(\pi \nu )-i\beta \cos (\pi \nu )}d\beta \right) dr \\
&&\overset{\nu \rightarrow 1}{\rightarrow }\frac{1}{2\pi }\int_{0}^{+\infty
}re^{-rt}\left( \int_{-\infty }^{+\infty }e^{-ir\beta +i\beta }d\beta
\right) dr \\
&=&\int_{0}^{+\infty }re^{-rt}\delta (r-1)dr=e^{-t}=E_{1,1}(-t).
\end{eqnarray*}

For large $t$, the following approximations follow from (\ref{rai8}) and (%
\ref{rai9}):%
\begin{equation*}
E_{\nu ,\beta }(-t^{\nu })\simeq \left\{
\begin{array}{l}
\frac{\Gamma (\nu -\beta +1)}{\pi t^{\nu }}\sin \left( (\beta -\nu )\pi
\right) ,\qquad \nu \neq \beta \\
\frac{\Gamma (\nu +1)}{\pi t^{2\nu }}\sin \left( \nu \pi \right) ,\qquad \nu
=\beta%
\end{array}%
\right. .
\end{equation*}

We study now a similar expansion for the Wright function, valid for any $%
0<\nu <1$ and $\beta \neq 1$. By following the same steps as in the proof of
(\ref{rai8}), we have that%
\begin{eqnarray*}
&&W_{\nu ,\beta }(-t^{\nu }) \\
&=&\sum_{m=0}^{\infty }\frac{(-t)^{\nu m}}{m!\Gamma (\nu m+\beta )} \\
&=&\frac{t^{1-\beta }}{2\pi i}\int_{0}^{\infty }dy\int_{0}^{\infty
}e^{-r(y+t)}y^{\beta -1}\left[ e^{i\pi \beta }W_{\nu ,\beta }(-y^{\nu
}e^{i\pi \nu })-e^{-i\pi \beta }W_{\nu ,\beta }(-y^{\nu }e^{-i\pi \nu })%
\right] dr \\
&=&\frac{t^{1-\beta }}{2\pi i}\left[ e^{i\pi \beta }\int_{0}^{\infty
}e^{-rt}\sum_{m=0}^{\infty }\frac{(-e^{i\pi \nu })^{m}}{m!\Gamma (\nu
m+\beta )}\frac{\Gamma (\nu m+\beta )}{r^{\beta +\nu m}}dr-e^{-i\pi \beta
}\int_{0}^{\infty }e^{-rt}\sum_{m=0}^{\infty }\frac{(-e^{-i\pi \nu })^{m}}{%
m!\Gamma (\nu m+\beta )}\frac{\Gamma (\nu m+\beta )}{r^{\beta +\nu m}}dr%
\right] \\
&=&\frac{t^{1-\beta }}{2\pi i}\int_{0}^{\infty }\frac{e^{-rt}}{r^{\beta }}%
\left[ e^{i\pi \beta -\frac{e^{i\pi \nu }}{r^{\nu }}}-e^{-i\pi \beta -\frac{%
e^{-i\pi \nu }}{r^{\nu }}}\right] dr \\
&=&\frac{t^{1-\beta }}{\pi }\int_{0}^{\infty }\frac{e^{-rt}}{r^{\beta }}e^{-%
\frac{\cos (\pi \nu )}{r^{\nu }}}\left[ \frac{e^{i\pi \beta -\frac{i\sin
(\pi \nu )\pi }{\pi r^{\nu }}}-e^{-i\pi \beta +\frac{i\sin (\pi \nu )\pi }{%
\pi r^{\nu }}}}{2i}\right] dr \\
&=&\frac{t^{1-\beta }}{\pi }\int_{0}^{\infty }\frac{e^{-rt}}{r^{\beta }}e^{-%
\frac{\cos (\pi \nu )}{r^{\nu }}}\sin \left( \pi \left( \beta -\frac{\sin
(\pi \nu )}{\pi r^{\nu }}\right) \right) dr.
\end{eqnarray*}

\

\end{document}